\documentclass[10pt]{article}

\usepackage{amsmath}
\usepackage{amssymb}
\usepackage{amstext}
\usepackage{amsthm}
\usepackage{a4wide}
\usepackage[top=2cm, bottom=2cm]{geometry}
\usepackage{color}
\usepackage{url}
\usepackage{graphicx}
\usepackage{algorithmic}
\usepackage{algorithm}
\usepackage{tikz}
\usepackage{placeins}
\usepackage{booktabs}

\definecolor{gray}{RGB}{150,150,150}
\definecolor{myblue}{RGB}{0,106,214}
\definecolor{unknown}{RGB}{0,100,200}

\DeclareMathOperator*{\argmin}{argmin}
\DeclareMathOperator{\myspan}{span}

\newtheorem{remark}{Remark}

\newcommand{\R} {\mathbb{R}}
\newcommand{\N} {\mathbb{N}}

\newcommand{\mcA} {\mathcal{A}}
\newcommand{\mcB} {\mathcal{B}}
\newcommand{\mcC} {\mathcal{C}}

\newcommand{\mcH} {\mathcal{H}}

\newcommand{\mcL} {\mathcal{L}}
\newcommand{\mcM} {\mathcal{M}}

\newcommand{\mcP} {\mathcal{P}}

\newcommand{\mcV} {\mathcal{V}}

\newcommand{\bF} {\mathbf{F}}
\newcommand{\bG} {\mathbf{G}}

\newcommand{\bS} {\mathbf{S}}

\newcommand{\bU} {\mathbf{U}}
\newcommand{\bV} {\mathbf{V}}

\newcommand{\bX} {\mathbf{X}}
\newcommand{\bY} {\mathbf{Y}}
\newcommand{\bZ} {\mathbf{Z}}
\newcommand{\br} {\mathbf{r}}

\newcommand{\dU} {\delta U}

\newcommand{\unl}{\mathsf{L}}

\newcommand{\T}{\mathsf{T}}
\newcommand{\minT}{{-\mathsf{T}}}

\DeclareMathOperator*{\tensor}{\raisebox{-.3ex}{\text{\huge$\times$}}}

\DeclareMathOperator{\tenvec}{vec}
\DeclareMathOperator{\grad}{grad}
\DeclareMathOperator{\rank}{rank}
\DeclareMathOperator{\rankT}{rank_{\text{ML}}}
\DeclareMathOperator{\rankTT}{rank_{\text{TT}}}
\DeclareMathOperator{\Proj}{P}

\newcommand{\Mr} {\mcM_\mathbf{r}}

\begin{document}

\title{Preconditioned low-rank Riemannian optimization for \\linear systems with tensor product structure}

\author{Daniel~Kressner\footnote{MATHICSE-ANCHP, \'Ecole Polytechnique F\'ed\'erale de Lausanne, Station 8, 1015 Lausanne, Switzerland. E-mail: daniel.kressner@epfl.ch}\and Michael~Steinlechner\footnote{MATHICSE-ANCHP, \'Ecole Polytechnique F\'ed\'erale de Lausanne, Station 8, 1015 Lausanne, Switzerland. E-mail: michael.steinlechner@epfl.ch
\newline
The work of M. Steinlechner has been supported by the SNSF research module \emph{Riemannian optimization for solving high-dimensional problems with low-rank tensor techniques} within the SNSF ProDoc \emph{Efficient Numerical Methods for Partial Differential Equations}. 
} \and~Bart~Vandereycken\footnote{Section de Math\'ematiques, Universit\'e de Gen\`eve, 2-4 rue du Li\`evre, 1211 Gen\`eve, Switzerland. E-Mail: bart.vandereycken@unige.ch}}

\date{July 28, 2015}
\maketitle

\begin{abstract}
The numerical solution of partial differential equations on high-dimensional domains gives
rise to computationally challenging linear systems. When using standard discretization techniques,
the size of the linear system grows exponentially with the number of dimensions, making the use of classic
iterative solvers infeasible. During the last few years, low-rank tensor approaches 
have been developed that allow to mitigate this \emph{curse of dimensionality} by exploiting the
underlying structure of the linear operator. In this work, we focus on tensors
represented in the Tucker and tensor train formats. We propose two preconditioned gradient methods on
the corresponding low-rank tensor manifolds: A Riemannian version of the
preconditioned Richardson method as well as an approximate Newton scheme
based on the Riemannian Hessian. For the latter, considerable attention is given
to the efficient solution of the resulting Newton equation.  In numerical
experiments, we compare the efficiency of our Riemannian algorithms with other
established tensor-based approaches such as a truncated preconditioned Richardson
method and the alternating linear scheme. The results show that
our approximate Riemannian Newton scheme is significantly faster in cases when the
application of the linear operator is expensive.


\end{abstract}

\section{Introduction}

This work is concerned with the approximate solution of large-scale linear systems
$Ax = f$ with $A \in \R^{n\times n}$. In certain applications, such as
the structured discretization of $d$-dimensional partial differential equations (PDEs),
the size of the linear system naturally decomposes as $n = n_1 n_2 \cdots n_d$ with $n_\mu \in \N$
for $\mu = 1,\ldots, d$. This allows us to view $Ax = f$ as a tensor equation
\begin{equation}\label{eq:calAxf}
    \mcA \bX = \bF,
\end{equation}
where $\bF, \bX \in \R^{n_1 \times n_2 \times \cdots \times n_d}$
are tensors of order $d$ and 
$\mcA$ is a linear operator on $\R^{n_1 \times n_2 \times \cdots \times n_d}$.

The tensor equations considered in this paper admit a decomposition of the form 
\begin{equation} \label{eq:laplace+potential}
 \mcA = \mcL + \mcV,
\end{equation}
where $\mcL$ is a Laplace-like operator with the matrix representation
\begin{equation} \label{eq:laplacelike}
  L = I_{n_d} \otimes \cdots \otimes I_{n_2} \otimes L_1 + 
 I_{n_d} \otimes \cdots \otimes I_{n_3} \otimes L_2  \otimes I_{n_1} + \cdots + 
 L_d \otimes I_{n_{d-1}} \otimes \cdots \otimes I_1,
\end{equation}
with matrices $L_\mu \in \R^{n_\mu \times n_\mu}$ and identity matrices $I_{n_\mu}$. 
The term $\mcV$ is dominated by $\mcL$ in the sense that $\mcL$ is assumed to be a good preconditioner 
for $\mcA$. Equations of this form arise, for example, from the discretization of the Schr\"odinger Hamiltonian~\cite{Lubich2008}, for which $\mcL$ and $\mcV$ correspond to the discretization of the kinetic and the potential energy terms, respectively. In this application, $\mcA$ (and thus also $L_\mu$) is symmetric positive definite. In the following, we restrict ourselves to this case, although some of the developments can, in principle, be generalized to indefinite and nonsymmetric matrices.

Assuming $\mcA$ to be symmetric positive definite allows us to reformulate~\eqref{eq:calAxf}
as an optimization problem 
\begin{equation}
    \label{eq:fullopt}
    \min_{{\bX} \in \R^{n_1 \times \cdots \times n_d}} \;\frac{1}{2} \langle \bX, \mcA \bX \rangle - \langle \bX, \bF \rangle
\end{equation}
It is well-known that the above problem is equivalent to minimizing the $\mcA$-induced norm of the error $\| \bX - \mcA^{-1} \bF \|_{\mcA}$. Neither~\eqref{eq:calAxf} nor~\eqref{eq:fullopt} are computationally tractable for larger values of $d$.
During the last decade, low-rank tensor techniques have been developed that aim at dealing with this curse of dimensionality by approximating ${\bf F}$ and ${\bf X}$ in a compressed format; see~\cite{Grasedyck2013,Hackbusch2012} for overviews. One approach consists of restricting~\eqref{eq:fullopt} to a subset $\mcM \subset \R^{n_1 \times n_2 \times \cdots \times n_d}$ of compressed tensors:
\begin{equation}
    \label{eq:problem_opt}
        \min_{{\bX} \in \mcM} \;f(\bX):=\frac{1}{2} \langle \bX, \mcA \bX \rangle - \langle \bX, \bF \rangle.
\end{equation}
Examples for $\mcM$ include the Tucker format \cite{Tucker1966,Kolda2009}, the tensor train (TT) format \cite{Oseledets2011a}, the matrix product states (MPS) format \cite{Affleck:1987} or the hierarchical Tucker format \cite{Grasedyck2010a,Hackbusch2009a}. Assuming that the corresponding ranks are fixed, $\mcM$ is a smooth embedded submanifold of $\R^{n_1 \times n_2 \times \cdots \times n_d}$ for each of these formats \cite{Holtz2010,Uschmajew2013,Uschmajew2013a,Haegeman:2014}. This property does not hold for the CP format, which we will therefore not consider.

When $\mcM$ is a manifold, Riemannian optimization techniques \cite{Absil2008} can be used to
address~\eqref{eq:problem_opt}. In a related context, 
first-order methods, such as Riemannian steepest descent and nonlinear CG,
have been successfully applied to matrix completion \cite{Boumal2011,Mishra2014a,Ngo2012,Vandereycken2012} and tensor completion \cite{DaSilva2015,Kressner2013a,Rauhut2015,Steinlechner2015}.

Similar to Euclidean optimization, the condition number of the Riemannian Hessian of the objective function is instrumental in predicting the performance of first-order optimization algorithms on manifolds; see, e.g.,~\cite[Thm.~2]{Luenberger:70kx} and~\cite[Thm.~4.5.6]{Absil2008}.
As will be evident from~\eqref{eq:Hessian} in \S\ref{sec:approx_Newton}, an ill-conditioned operator $\mcA$ can be expected to yield an ill-conditioned Riemannian Hessian. As this is the case for the applications we consider, any naive first-order method will be prohibitively slow and noncompetitive with existing methods. 

For Euclidean optimization, it is well known that preconditioning or,
equivalently, adapting the underlying metric can be used to address the slow
convergence of such first-order methods.  Combining steepest descent with the
Hessian as a (variable) preconditioner yields the Newton method with (local)
second order convergence~\cite[Sec. 1.3.1]{Nesterov2004}. To overcome the high
computational cost associated with Newton's method, several approximate Newton
methods exist that use  cheaper second-order models. For example, Gauss--Newton
is a particularly popular approximation when solving non-linear least-squares
problems. For Riemannian optimization, the connection between preconditioning
and adapting the metric is less immediate and we  explore both directions to
speed up first-order methods. On the one hand, we will consider a rather ad hoc
way to precondition the Riemannian gradient direction. On the other hand, we
will consider an approximate Newton method that can be interpreted as a
constrained Gauss--Newton method. This requires setting up and solving linear
systems with the Riemannian Hessian or an approximation thereof.
In~\cite{Vandereycken2010}, it was shown that neglecting curvature terms in the
Riemannian Hessian leads to an efficient low-rank solver for Lyapunov matrix
equations. We will extend these developments to more general equations with
tensors approximated in the Tucker and the TT formats.

Riemannian optimization is by no means the only sensible approach to finding
low-rank tensor approximations to the solution of the linear
system~\eqref{eq:calAxf}. For linear operators only involving the Laplace-like
operator~\eqref{eq:laplacelike}, exponential sum
approximations~\cite{Grasedyck2004,Hackbusch2006} and tensorized Krylov
subspace methods~\cite{Kressner2010} are effective and allow for a thorough
convergence analysis.  For more general equations, a straightforward approach
is to apply standard iterative methods, such as the Richardson iteration or the
CG method, to~\eqref{eq:calAxf} and represent all iterates in the low-rank
tensor format;
see~\cite{Ballani2013,Dolgov2012h,Khoromskij2010b,Khoromskij2011c,Kressner2011}
for examples. One critical issue in this approach is to strike a balance
between maintaining convergence and avoiding excessive intermediate rank growth
of the iterates. Only recently, this has been analyzed in more detail
\cite{Bachmayr2014}.  A very different approach consists of applying
alternating optimization techniques to the constrained optimization
problem~\eqref{eq:problem_opt}. Such methods have originated in quantum
physics, most notably the so called
DMRG method to address eigenvalue problems in the context of strongly
correlated quantum lattice systems, see~\cite{Schollwock2011} for an overview.
The ideas of DMRG and related methods have been extended to linear systems in
the numerical analysis community
in~\cite{Dolgov2012c,Dolgov2014a,Holtz2012,Oseledets2011c} and are generally
referred to as alternating linear schemes (ALS).  While such methods often
exhibit fast convergence, especially for operators of the
form~\eqref{eq:laplace+potential}, their global convergence properties are
poorly understood. Even the existing local convergence results for
ALS~\cite{Rohwedder2013,Uschmajew2013a} offer little intuition on the
convergence \emph{rate}.  The efficient implementation of ALS for low-rank
tensor formats can be a challenge. In the presence of larger ranks, the (dense)
subproblems that need to be solved in every step of ALS are large and tend to
be ill-conditioned. In~\cite{Kressner2013b,Kressner2011a}, this issue has been
addressed by combining an iterative solver with a preconditioner tailored to
the subproblem. The design of such a preconditioner is by no means simple, even
the knowledge of an effective preconditioner for the full-space
problem~\eqref{eq:calAxf} is generally not sufficient.  So far, the only known
effective preconditioners are based on exponential sum approximations for
operators with Laplace-like structure~\eqref{eq:laplacelike}, which is
inherited by the subproblems.

Compared to existing approaches, the preconditioned low-rank Riemannian optimization
methods proposed in this paper have a number of advantages. Due to imposing the manifold
constraint, the issue of rank growth is completely avoided. Our methods have a global nature,
all components of the low-rank tensor format are improved at once and hence the stagnation
typically observed during ALS sweeps is avoided. Moreover, we completely  avoid the need
for solving subproblems very accurately. One of our methods can
make use of preconditioners for the full-space problem~\eqref{eq:calAxf}, while for the other methods
preconditioners are implicitly obtained from approximating the Riemannian Hessian.
A disadvantage shared with existing methods, our method strongly relies on the 
decomposition~\eqref{eq:laplace+potential} of the operator to construct effective preconditioners.

In passing, we mention that there is another notion of preconditioning for
Riemannian optimization on low-rank matrix manifold, see,
e.g.,~\cite{Mishra2014b,Mishra2014,Ngo2012}.  These techniques address the
ill-conditioning of the manifold parametrization, an aspect that is not related
and relevant to our developments, as we do not directly work with the
parametrization.

The rest of this paper is structured as follows. In Section~\ref{sec:manifold},
we briefly review the Tucker and TT tensor formats and the structure of the
corresponding manifolds. Then, in Section~\ref{sec:first-order}, a Riemannian
variant of the preconditioned Richardson method is introduced. In
Section~\ref{sec:second-order}, we incorporate second-order information using
an approximation of the Riemannian Hessian of the cost function and solving the
corresponding Newton equation.  Finally, numerical experiments comparing the
proposed algorithms with existing approaches are presented in
Section~\ref{sec:experiments}.

\section{Manifolds of low-rank tensors}
\label{sec:manifold}

In this section, we discuss two different representations for tensors $\bX \in \R^{n_1 \times
    n_2 \times \cdots \times n_d}$, namely the \emph{Tucker} and \emph{tensor train/matrix
    product states} (TT/MPS)
formats, along with their associated notions of low-rank structure and their geometry.  We
will only mention the main results here and refer to the articles by Kolda and Bader
\cite{Kolda2009} and by Oseledets
\cite{Oseledets2011a} for more details.
More elaborate discussions on the manifold structure and computational
efficiency considerations can be found in \cite{Koch2010,Kressner2013a} for the
Tucker format and in \cite{Lubich2014,Steinlechner2015,Uschmajew2013a} for the
TT format, respectively.

\subsection{Tucker format}\label{sec:Tucker_manifold}

\paragraph{Format.} The \emph{multilinear rank} of a tensor $\bX \in \R^{n_1 \times n_2 \times \cdots \times n_d}$ is
defined as the $d$-tuple
\[
    \rankT(\bX) = (r_1, r_2,\ldots, r_d) = \left(\rank(\bX_{(1)}),
        \rank(\bX_{(2)}),\ldots,\rank(\bX_{(d)})\right)
\]
with 
\[
 \bX_{(\mu)} \in \R^{n_\mu \times (n_1 \cdots n_{\mu-1} n_{\mu+1} \cdots n_d)}, \qquad
 \mu = 1,\ldots, d,
\]
the $\mu$th matricization of $\bX$; see \cite{Kolda2009} for more details. 

Any tensor $\bX
\in \R^{n_1 \times n_2 \times \cdots \times n_d}$ of multilinear rank $\br = (r_1, r_2,\ldots, r_d) $ can be represented as 
\begin{equation} \label{eq:tuckerformat}
     \bX( i_1, \ldots, i_d ) = \sum_{j_1=1}^{r_1} \cdots \sum_{j_{d-1}=1}^{r_{d-1}}
    \bS(j_1,j_2,\ldots,j_d) U_1(i_1,j_1) U_2(i_2, j_3) \cdots U_d(i_{d-1},j_d),
\end{equation}
for some \emph{core tensor} $\bS \in \R^{r_1 \times \cdots \times r_d}$ and \emph{factor
matrices} $U_\mu \in \R^{n_\mu \times r_\mu}, \mu = 1,\ldots,d$. In the following, we always choose the factor matrices to have orthonormal columns, $U_\mu^\T U_\mu = I_{r_\mu}$.

Using the $\mu$th mode product $\times_\mu$, see~\cite{Kolda2009}, 
one can write~\eqref{eq:tuckerformat} more compactly as
\begin{equation} \label{eq:tuckerformat2}
    \bX = \bS \times_1 U_1 \times_2 U_2\  \cdots \times_d U_d.
\end{equation}

\paragraph{Manifold structure.} It is known
\cite{Koch2010,Hackbusch2012,Uschmajew2013} that the set of tensors having
multilinear rank $\br$ forms a smooth submanifold embedded in $\R^{n_1 \times
    n_2 \times \cdots \times n_d}$. This manifold $\Mr$ is of dimension
\[
    \dim \Mr =  \prod_{\mu=1}^{d-1} r_\mu + \sum_{\mu=1}^d r_{\mu} n_\mu - r_{\mu}^2.
\]
For $\bX \in \Mr$ represented as in~\eqref{eq:tuckerformat2},
any tangent vector $\xi \in T_{\bX} \Mr$ can be written as
\begin{equation} \label{eq:tangenttuckerrepresentation}
\begin{aligned}
  \xi = \bS \times_1 &\delta U_1 \times_2 U_2\  \cdots \times_d U_d 
  \;+\; \bS \times_1 U_1 \times_2 \delta U_2\  \cdots \times_d U_d\\
  &\;+\; \cdots 
  \;+\; \bS \times_1 U_1 \times_2 U_2\  \cdots \times_d \delta U_d
  \;+\; \delta\bS \times_1 U_1 \times_2 U_2\  \cdots \times_d  U_d,
\end{aligned}
\end{equation}
for some first-order variations $\delta \bS \in \R^{r_1 \times \cdots \times r_d}$ and
$\delta U_\mu \in \R^{n_\mu \times r_\mu}$.
This representation of tangent vectors allows us to decompose the tangent
space $T_{\bX} \Mr$ orthogonally as 
\begin{equation} \label{eq:tangentdecomptucker}
     T_\bX \mcM = \mcV_1 \oplus \mcV_2 \oplus \cdots \oplus \mcV_d \oplus \mcV_{d+1}, \quad \text{with } \mcV_\mu \perp \mcV_\nu \;\, \forall\, \mu \neq \nu,
\end{equation}
where the subspaces $\mcV_\mu$ are given by
\begin{equation}\label{eq:gauged_V_Tucker}
    \mcV_\mu = \Big\{   \bS \times_\mu
      \delta U_\mu \tensor_{\substack{\nu=1 \\ \nu \neq \mu}}^d  U_\nu
	\colon \delta U_\mu \in \R^{n_\mu \times r_\mu}, \ \delta U_\mu^\T  U_\mu = 0
	\Big\}, \qquad \mu = 1,\ldots,d,
\end{equation}
and 
\[
    \mcV_{d+1} = \left\{
	\delta\bS \tensor_{\nu=1}^d  U_\nu \colon \delta \bS \in \R^{r_1 \times \cdots \times r_d}
                \right\}.
\]
In particular, this decomposition shows that, given the core tensor $\bS$ and
factor matrices $U_\mu$ of $\bX$, the tangent vector $\xi$ is uniquely
represented in terms of $\delta \bS$  and \emph{gauged} $\delta U_\mu$ .

\paragraph{Projection onto tangent space.}
Given $\bZ \in \R^{n_1 \times \cdots \times n_d}$, the components $\delta
U_\mu$ and $\delta \bS$ of the orthogonal projection $\xi =
\Proj_{T_{\bX}\Mr}(\bZ)$  are given by (see \cite[Eq.(2.7)]{Koch2010}) 
\begin{equation}\label{eq:projectionTucker}
	\begin{aligned}
 \delta \bS &= \bZ \tensor_{\mu=1}^d U_\mu^\T,\\
\delta U_\mu &= (I_{n_\mu} - U_\mu U_\mu^\T) \Big[ \bZ \tensor_{\substack{\nu=1 \\ \nu \neq \mu}}^d U_\nu^\T \Big]_{(1)} \bS_{(\mu)}^\dagger,
	\end{aligned}
\end{equation}
where $\bS_{(\mu)}^\dagger = \bS_{(\mu)}^\T \big(
\bS_{(\mu)}\bS_{(\mu)}^\T\big)^{-1}$ is the Moore--Penrose pseudo-inverse of
$\bS_{(\mu)}$. 
The projection of a Tucker tensor of multilinear rank $\tilde{\mathbf{r}}$ into $T_\bX \Mr$ can
be performed in $O( dn \tilde{r} r^{d-1} + \tilde{r}^{d} r)$ operations, where we set $\tilde{r} := \max_\mu
\tilde{r}_\mu$, $r := \max_\mu r_\mu$ and $\tilde{r} \geq r$.  

\subsection{Representation in the TT format} \label{sec:intro_TT}
\paragraph{Format.} The TT format is (implicitly) based on matricizations that
merge the first $\mu$ modes into row indices and the remaining indices into column indices:
\[
 \bX^{<\mu>} \in \R^{(n_1\cdots n_\mu)\times (n_{\mu+1} \cdots n_d)}, \qquad
 \mu = 1,\ldots, d-1.
\]
The \emph{TT rank} of $\bX$ is the tuple $\rankTT(\bX) := 
(r_0, r_1, \ldots, r_d)$ with $r_\mu = \rank(\bX^{<\mu>})$.
By definition, $r_0 = r_d = 1$ .

A tensor $\bX\in \R^{n_1 \times n_2 \times \dots \times n_d}$ of TT rank 
$\br = (r_0, r_1, \ldots, r_d)$ admits the representation 
\begin{equation}
    \label{eq:TT_pointwise}
    \bX( i_1, \ldots, i_d ) = U_1(i_1) U_2(i_2) \cdots U_d(i_d)
\end{equation}
where each $U_\mu(i_\mu)$ is  a matrix of size $r_{\mu-1} \times r_{\mu}$
for $i_\mu = 1,2,\dots,n_\mu$. By stacking the matrices
$U_\mu(i_\mu)$, $i_\mu = 1,2,\dots,n_\mu$ into third-order tensors 
$\bU_\mu$ of size $r_{\mu-1} \times n_\mu \times r_\mu$, the so-called 
\emph{TT cores}, we can also write 
\eqref{eq:TT_pointwise} as
\[	
    \bX( i_1, \ldots, i_d ) = \sum_{j_1=1}^{r_1} \cdots \sum_{j_{d-1}=1}^{r_{d-1}}
    \bU_1(1,i_1,j_1) \bU_2(j_2, i_2, j_3) \cdots \bU_d(j_{d-1},i_d,1).
\]
To access and manipulate individual cores, it is useful to introduce the \emph{interface matrices}
\begin{align*}
	\bX_{\le \mu} &= [U_1(i_1)U_2(i_2) \cdots U_\mu(i_\mu)] \in \R^{n_1n_2\cdots n_\mu \times r_\mu},\\
	\bX_{\ge \mu} &= [U_\mu(i_{\mu})U_{\mu+1}(i_{\mu+1}) \cdots U_d(i_d)]^\T \in \R^{n_{\mu} n_{\mu+1} \cdots n_d \times r_{\mu-1}},
\end{align*}
and\begin{equation}\label{eq:Xneqmu}
\bX_{\neq \mu} = \bX_{\ge {\mu+1}} \otimes I_{n_\mu} \otimes \bX_{\le \mu-1}
\in \R^{n_1 n_2\cdots n_d \times r_{\mu-1} n_\mu r_{\mu}}.
\end{equation}
In particular, this allows us to pull out the $\mu$th core as $\tenvec(\bX) =
\bX_{\neq \mu} \tenvec(\bU_\mu)$, where $\tenvec(\cdot)$ denotes the
vectorization of a tensor.

There is some freedom in choosing the cores in the representation
\eqref{eq:TT_pointwise}. In particular, we can orthogonalize parts of $\bX$. We
say that $\bX$ is \emph{$\mu$-orthogonal} if $\bX_{\leq \nu}^\T \bX_{\leq \nu}
= I_{r_\nu}$ for all $\nu = 1,\ldots,\mu-1$ and $\bX_{\geq \nu} \bX_{\geq
    \nu}^\T = I_{r_{\nu-1}}$ for all $\nu = \mu+1,\ldots,d$, see,
e.g.,~\cite{Steinlechner2015} for more details.
\paragraph{Manifold structure.}
The set of tensors having fixed TT rank,
\[
    \Mr = \left\{ \bX \in \R^{n_1 \times \cdots \times n_d} \colon \rankTT(\bX) =
        \mathbf{r} \right\},
\]
forms a smooth embedded submanifold of $\R^{n_1 \times \cdots \times n_d}$,
see \cite{Holtz2010,Hackbusch2012,Uschmajew2013a}, 
of dimension
\[
    \dim \Mr = \sum_{\mu=1}^d r_{\mu-1}n_\mu r_{\mu} - \sum_{\mu=1}^{d-1} r_\mu^2.
\]
Similar to the Tucker format, the tangent space $T_\bX \Mr$ at $\bX \in \Mr$ admits an orthogonal decomposition:
\begin{equation} \label{eq:tangentdecomp}
     T_\bX \Mr = \mcV_1 \oplus \mcV_2 \oplus \cdots \oplus \mcV_d, \quad \text{with } \mcV_\mu \perp \mcV_\nu \;\, \forall\, \mu \neq \nu.
\end{equation}
Assuming that $\bX$ is $d$-orthogonal, the subspaces $\mcV_\mu$ can be represented as
\begin{equation}\label{eq:V_TT}
	\begin{aligned}
    \mcV_\mu &= \left\{ \,\bX_{\neq \mu} \tenvec(\delta \bU_\mu)\colon \delta \bU_\mu \in \R^{r_{\mu-1} \times n_\mu \times r_\mu}, 
                \big(\bU_\mu^{\unl}\big)^\T \delta \bU_\mu^\unl = 0  \right\}, \quad  \mu = 1,\ldots,d-1, \\				
    \mcV_d &= \left\{ \,\bX_{\neq d} \tenvec(\delta \bU_d)\colon \delta \bU_d \in \R^{r_{d-1} \times n_d \times r_d} 
                \right\}.
\end{aligned}
\end{equation}
Here, $\bU_\mu^{\unl} \equiv \bU_\mu^{<2>} \in \R^{r_{\mu-1} n_\mu \times r_\mu}$ is called the
\emph{left unfolding} of $\bU_\mu$ and it has orthonormal columns for $\mu = 1,\ldots,d-1$, due
to the $d$-orthogonality of $\bX$. The \emph{gauge conditions} $(\bU_\mu^{\unl})^\T \delta
\bU_\mu^\unl = 0$ for $\mu\neq d$ ensure the mutual orthogonality of the subspaces $\mcV_\mu$
and thus yield a unique representation of a tangent vector $\xi$ in terms of gauged $\delta
\bU_\mu$. Hence, we can write any tangent vector $\xi \in T_\bX \Mr$ in the convenient form
\begin{equation}
    \label{eq:tangentTTrepresentation}
    \xi = \sum_{\mu=1}^d \bX_{\neq \mu} \tenvec(\delta \bU_\mu) \;\in \R^{n_1 n_2 \cdots n_d} \qquad \text{s.t.} \qquad (\bU_\mu^{\unl})^\T \delta \bU_\mu^\unl = 0, \quad \forall \mu \neq d.
\end{equation}

\paragraph{Projection onto tangent space.} 

The orthogonal projection $\Proj_{T_\bX \mcM}$ onto the tangent space $T_\bX \mcM$ can be decomposed in accordance with~\eqref{eq:tangentdecomp}: 
\begin{equation*}
    \label{eq:projtangentdecomp}
     \Proj_{T_\bX \mcM} = \Proj^1 + \Proj^2 + \cdots + \Proj^d,
\end{equation*}
where $\Proj^\mu$ are orthogonal projections onto $\mcV_\mu$. Let $\bX\in \Mr$ be $d$-orthogonal and $\bZ \in \R^{n_1 \times \cdots \times n_d}$. Then the projection can be written as
\begin{equation}\label{eq:P_Z_TT_in_components}
 \Proj_{T_{\bX}\Mr}(\bZ)  = \sum_{\mu=1}^d \Proj^\mu (\bZ) \quad \text{where} \quad   \Proj^\mu (\bZ) = \bX_{\neq \mu} \tenvec(\delta \bU_\mu).
\end{equation}
For $\mu=1,\ldots,d-1$, the components $\delta \bU_\mu$ in this expression are given by \cite[p.~924]{Lubich2015}
\begin{equation}
    \label{eq:deltaU}
    \delta \bU_\mu^\unl = (I_{n_\mu r_{\mu-1}} - \Proj^\unl_\mu)\big(I_{n_\mu} \otimes \bX_{\leq \mu-1}^\T\big) \bZ^{<\mu>} \bX_{\geq \mu+1}\big( \bX_{\geq \mu+1}^\T \bX_{\geq \mu + 1}\big)^{-1}  
\end{equation}
with $ \Proj_\mu^\unl = \bU_\mu^\unl (\bU_\mu^\unl)^\T$ the orthogonal projector onto the range of $\bU_\mu^\unl$. For $\mu=d$, we have
\begin{equation}
    \label{eq:deltaU_d}
    \delta \bU_d^\unl = \big(I_{n_d} \otimes \bX_{\leq d-1}^\T\big) \bZ^{<d>}. 
\end{equation}
The projection of a tensor of TT rank $\tilde{\mathbf{r}}$ into $T_\bX \Mr$ can
be performed in $O( d n r \tilde{r}^2)$ operations, where we set again $\tilde{r} := \max_\mu
\tilde{r}_\mu$, $r := \max_\mu r_\mu$ and $\tilde{r} \geq r$.  

\begin{remark}\label{rem:inv_gram}
    Equation \eqref{eq:deltaU} is not well-suited for numerical calculations due to the presence of
the inverse of the Gram matrix $\bX_{\geq \mu+1}^\T \bX_{\geq \mu + 1}$, which is typically severely ill-conditioned. In
\cite{Khoromskij2012,Steinlechner2015}, it was shown that by $\mu$-orthogonalizing the $\mu$th summand of the
tangent vector representation, these inverses can be avoided at no extra
costs. To keep the notation short, we do not include this individual orthogonalization in the
equations above, but make use of it in the implementation of the algorithm and the numerical
experiments discussed in Section \ref{sec:experiments}.
\end{remark}

\subsection{Retractions} \label{sec:retraction}

Riemannian optimization algorithms produce search directions that are contained in
the tangent space $T_\bX \Mr$ of the current iterate. To obtain the next iterate
on the manifold, tangent vectors are mapped back to the manifold by application of a
\emph{retraction} map $R$ that satisfies certain properties; see~\cite[Def. 1]{Absil2012a} for a formal definition.

It has been shown in~\cite{Kressner2013a} that the higher-order SVD
(HOSVD)~\cite{DeLathauwer2000}, which aims at approximating a given tensor of rank $\tilde{\br}$ by a tensor of lower 
rank $\br$, constitutes a retraction on the Tucker manifold $\Mr$ that can be computed
efficiently in $O(dn\tilde{r}^2 + \tilde{r}^{d+1})$ operations. For the TT manifold, we will use the analogous
TT-SVD \cite[Sec. 3]{Oseledets2011a} for a retraction with a computational cost of $O(dn\tilde{}^3)$,
see~\cite{Steinlechner2015}. For both manifolds, we will denote by $R\big(\bX + \xi \big)$ the
retraction\footnote{Note that the domain of definition of $R$ is the affine
tangent space $\bX + T_\bX \Mr$, which departs from the usual convention to
define $R$ on $T_\bX \Mr$ and only makes sense for this particular type of
retraction.} of $\xi \in T_\bX \Mr$ that is computed by the HOSVD/TT-SVD of
$\bX + \xi$.

\section{First-order Riemannian optimization and preconditioning}
\label{sec:first-order}

In this section, we discuss ways to incorporate preconditioners into simple first-order Riemannian optimization methods.

\subsection{Riemannian gradient descent}
\label{subsec:RiemannianGD}

To derive a first-order optimization method on a manifold $\Mr$, we first need to construct the 
Riemannian gradient.
For the cost function~\eqref{eq:problem_opt} associated with linear systems, the Euclidean gradient is given by 
\[
    \nabla f(\bX) = \mcA \bX - \bF.
\]
For both the Tucker and the TT formats, $\Mr$ is an embedded submanifold of $\R^{n_1\times \cdots \times
    n_d}$ and hence the Riemannian gradient can be obtained by projecting $\nabla f$ onto the tangent space:
\[
    \grad f(\bX) = \Proj_{T_\bX \Mr} ( \mcA \bX - \bF).
\]
Together with the retraction $R$ of Section~\ref{sec:retraction}, this
yields the basic Riemannian gradient descent algorithm:
\begin{equation}
    \label{eq:riemannsteep}
    \bX_{k+1} = R\big(\bX_{k} + \alpha_k  \xi_k\big), \quad \text{with} \quad  
    \xi_k = -\Proj_{T_{\bX_k} \mcM} \nabla f(\bX_k).
\end{equation}
As usual, a suitable step size $\alpha_k$ is obtained by 
standard Armijo backtracking linesearch. Following~\cite{Vandereycken2012}, a good initial guess for the backtracking procedure is constructed by
an exact linearized linesearch on the tangent space alone (that is, by neglecting the retraction):
\begin{equation}
    \label{eq:linearizedlinesearch}
    \argmin_\alpha f(\bX_k + \alpha \xi_k ) = -\frac{\langle \xi_k, \nabla f(\bX_k)\rangle}{\langle \xi, \mcA \xi\rangle}. 
\end{equation}

\subsection{Truncated preconditioned Richardson iteration}

\begin{paragraph}{Truncated Richardson iteration.}
The Riemannian gradient descent defined by \eqref{eq:riemannsteep} closely resembles a
truncated Richardson iteration for solving linear systems:
\begin{equation}
    \label{eq:richardson}
    \bX_{k+1} = R\big(\bX_{k} + \alpha_k \xi_k \big), \quad \text{with} \quad  
    \xi_k = - \nabla f(\bX_k) = \bF - \mcA \bX_k,
\end{equation}
which was proposed for the CP tensor format in~\cite{Khoromskij2011c}.  For the
hierarchical Tucker format, a variant of the TT format, the
iteration~\eqref{eq:richardson} has been analyzed in~\cite{Bachmayr2014}.  In
contrast to manifold optimization, the rank does not need to be fixed but can
be adjusted to strike a balance between low rank and convergence speed. It has
been observed, for example in~\cite{Kressner2012}, that such an
iterate-and-truncate strategy greatly benefits from preconditioners, not only
to attain an acceptable convergence speed but also to avoid excessive rank
growth of the intermediate iterates. 
\end{paragraph}

\begin{paragraph}{Preconditioned Richardson iteration.}
For the standard Richardson iteration 
$\bX_{k+1} = \bX_{k} - \alpha_k \xi_k$, a symmetric positive definite preconditioner $\mcP$ 
for $\mcA$ can be incorporated as follows:
\begin{equation}\label{eq:prec_Richardson}
    \bX_{k+1} = \bX_{k} + \alpha_k \mcP^{-1} \xi_k \quad \text{with} \quad
    \xi_k = \bF - \mcA \bX_k.
\end{equation}
Using the Cholesky factorization $\mcP = \mcC \mcC^\T$, this iteration turns out to be equivalent to applying the Richardson iteration to the transformed
symmetric positive definite linear system
\begin{equation*}\label{eq:kacksystem}
 \mcC^{-1} \mcA \mcC^\minT \bY = \mcC^{-1} \bF
\end{equation*}
after changing coordinates by $\mcC^\T \bX_k$. At the same time,~\eqref{eq:prec_Richardson} can be viewed as applying gradient descent in the inner product induced by $\mcP$.
\end{paragraph}

\begin{paragraph}{Truncated preconditioned Richardson iteration.}
The most natural way of combining truncation with preconditioning leads to the \emph{truncated preconditioned Richardson iteration}
\begin{equation}\label{eq:trunc_prec_Richardson}
    \bX_{k+1} = R\big(\bX_{k} + \alpha_k \mcP^{-1} \xi_k \big), \quad \text{with} \quad 
    \xi_k = \bF - \mcA \bX_k,
\end{equation}
see also~\cite{Khoromskij2011c}.
In view of Riemannian gradient descent~\eqref{eq:riemannsteep}, it appears natural to
project the search direction to the tangent space,
leading to the ``geometric'' variant
\begin{equation}\label{eq:trunc_prec_projected_Richardson}
    \bX_{k+1} = R\big(\bX_{k} + \alpha_k \Proj_{T_{\bX_k} \Mr} \mcP^{-1} \xi_k \big), \quad \text{with}  \quad 
    \xi_k = \bF - \mcA \bX_k.
\end{equation}

In terms of convergence, we have observed that the
scheme~\eqref{eq:trunc_prec_projected_Richardson} behaves similar
to~\eqref{eq:trunc_prec_Richardson}; see \S\ref{sec:experiments_tucker}.
However, it can be considerably cheaper per iteration: Since only tangent
vectors need to be retracted in \eqref{eq:trunc_prec_projected_Richardson}, the
computation of the HOSVD/TT-SVD in $R$ involves only tensors of bounded rank,
regardless of the rank of $\mcP^{-1} \xi_k$. In particular, with $\br$ the
Tucker/TT rank of $\bX_k$, the corresponding rank of $\bX_{k} + \alpha_k
\Proj_{T_{\bX_k} \Mr} \mcP^{-1} \xi_k$ is at most $2 \br$; see
\cite[\S3.3]{Kressner2013a} and \cite[Prop.~3.1]{Steinlechner2015}. On the
other hand, in \eqref{eq:trunc_prec_Richardson} the rank of $\bX_{k} + \alpha_k
\mcP^{-1} \xi_k$ is determined primarily by the quality of the preconditioner
$\mcP$ and can possibly be very large. 

Another advantage occurs for the special but important case when $\mcP^{-1} =
\sum_{\alpha=1}^s\mcP_\alpha$, where each term $\mcP_\alpha$ is relatively
cheap to apply. For example, when $\mcP^{-1}$ is an exponential sum
preconditioner~\cite{Braess2005} then $s = d$ and $\mcP_\alpha$ is a Kronecker
product of small matrices. By the linearity of $\Proj_{T_{\bX_k} \Mr}$, we have
\begin{equation}
    \label{eq:proj_commute_prec}
 \Proj_{T_{\bX_k} \Mr} \mcP^{-1} \xi_k = \sum_{\alpha=1}^s \Proj_{T_{\bX_k} \Mr} \mcP_\alpha \xi_k,
\end{equation}
which makes it often cheaper to evaluate this expression in the iteration~\eqref{eq:trunc_prec_projected_Richardson}.
To see this, for example, for the TT
format, suppose that $\mcP_\alpha \xi$ has TT ranks $r_p$. Then the preconditioned direction $\mcP^{-1} \xi_k$ can be expected
to have TT ranks $s r_p$.
Hence, the straightforward application of $\Proj_{T_{\bX_k} \Mr}$ to $\mcP^{-1} \xi_k$
requires $O(dn (s r_p)^2 r)$ operations. 
Using the expression on the right-hand side of~\eqref{eq:proj_commute_prec} instead reduces the cost to
$O(dn s  r_p^2 r)$ operations, since the summation of tangent vectors amounts to simply adding
their parametrizations.
In contrast, since the retraction is a non-linear operation,
trying to achieve similar cost savings in~\eqref{eq:trunc_prec_Richardson} by simply truncating the
culmulated sum subsequently may lead to
severe cancellation \cite[\S6.3]{Kressner:2014}.

\end{paragraph}

\section{Riemannian optimization using a quadratic model}
\label{sec:second-order}

As we will see in the numerical experiments in Section~\ref{sec:experiments},
the convergence of the first-order methods presented above crucially depends on the availability
of a good preconditioner for the full problem. 
In this section, we present Riemannian optimization methods based on a quadratic model. In these
methods, the preconditioners are derived from an approximation of the Riemannian Hessian.

\subsection{Approximate Newton method}\label{sec:approx_Newton}

The Riemannian Newton method~\cite{Absil2008} applied to~\eqref{eq:problem_opt}
determines the search direction $\xi_k$ from the equation
\begin{equation}\label{eq:Newton}
 \mcH_{\bX_k} \xi_k = - \Proj_{T_\bX \Mr} \nabla f(\bX_k),
\end{equation}
where the symmetric linear operator $\mcH_{\bX_k}:T_{\bX_k} \Mr \to T_{\bX_k} \Mr$ is the Riemannian Hessian of \eqref{eq:problem_opt}. Using  \cite{AbsMahTru2013}, we have
\begin{align}
 \mcH_{\bX_k} &= \Proj_{T_{\bX_k} \Mr} \big[ \nabla^2 f(\bX_k) +  J_{\bX_k} \nabla^2 f(\bX_k) \big]  \Proj_{T_{\bX_k} \Mr} \notag\\
 &= \Proj_{T_{\bX_k} \Mr} \big[ \mcA +  J_{\bX_k} (\mcA \bX_k - \bF) \big]  \Proj_{T_{\bX_k} \Mr} \label{eq:Hessian}
\end{align}
with the Fr\'echet derivative\footnote{$J_{\bX_k}$ is an operator from $\R^{n
        \times n \times \cdots \times n}$ to the space of self-adjoint linear
    operators $T_{\bX_k} \Mr \to T_{\bX_k} \Mr$.} $J_{\bX_k}$  of
$\Proj_{T_{\bX_k} \Mr}$.

As usual, the Newton equation is only well-defined near a strict local
minimizer and solving it exactly is prohibitively expensive in a large-scale
setting. We therefore approximate the linear system~\eqref{eq:Newton} in two
steps: First, we drop the term containing $J_{\bX_k}$ and second, we replace
$\mcA = \mcL + \mcV$ by $\mcL$. The first approximation can be interpreted as
neglecting the curvature of $\Mr$, or equivalently, as linearizing the manifold
at $\bX_k$. Indeed, this term is void if $\Mr$ would be a (flat) linear
subspace. This approximation is also known as constrained Gauss--Newton (see,
e.g,~\cite{Bock_1987}) since it replaces the constraint $\bX \in \Mr$ with its
linearization $\bX \in T_\bX \Mr$ and neglects the constraints in the
Lagrangian. The second approximation is natural given the assumption of $\mcL$
being a good preconditioner for $\mcA = \mcL + \mcV$. In addition, our
derivations and numerical implementation will rely extensively on the fact that
the Laplacian $\mcL$ acts on each tensor dimension separately.

The result is an \emph{approximate Newton method} were the search direction $\xi_k$ is determined from
\begin{equation}\label{eq:GN_manifold}
 \Proj_{T_{\bX_k} \Mr} \mcL \Proj_{T_{\bX_k} \Mr} \xi_k = \Proj_{T_\bX \Mr} ( \bF - \mcA \bX_k).
\end{equation}
Since $\mcL$ is positive definite, this equation is always well-defined for any $\bX_k$. In addition, $\xi_k$ is also gradient-related and hence the iteration
\begin{equation*}
    \bX_{k+1} = R\big(\bX_{k} + \alpha_k \xi_k \big)
\end{equation*}
is guaranteed to converge globally to a stationary point of the cost function if $\alpha_k$ is determined from Armijo backtracking \cite{Absil2008}.

Despite all the simplifications, the numerical solution
of~\eqref{eq:GN_manifold} turns out to be a nontrivial task. In the following
section, we explain an efficient algorithm for solving \eqref{eq:GN_manifold}
exactly when $\Mr$ is the Tucker manifold. For the TT manifold, this approach
is no longer feasible and we therefore present an effective preconditioner that
can used for solving~\eqref{eq:GN_manifold} with the preconditioned CG method.

\subsection{The approximate Riemannian Hessian in the Tucker case}
\label{sec:newton_tucker}
The solution of 
the linear system~\eqref{eq:GN_manifold} 
was addressed for the matrix case ($d=2$) in~\cite[Sec. 7.2]{Vandereycken2010}.
In the following, we extend this approach to tensors in the Tucker format.
To keep the presentation concise, we restrict ourselves to $d=3$; the extension to $d>3$ is straightforward.

For tensors of order $3$ in the Tucker format, we write~\eqref{eq:GN_manifold} as follows:
\begin{equation}\label{eq:GN_tucker}
    \Proj_{T_\bX \Mr} \mcL \Proj_{T_\bX \Mr} \textcolor{unknown}{\xi} =  \eta,
\end{equation}
where
\begin{itemize}
 \item $\bX \in \Mr$ is parametrized by factor matrices $U_1,U_2,U_3$ having \emph{orthonormal columns} and the core tensor $\bS$;
 \item the right-hand side $\eta \in {T_\bX \Mr}$ is given in terms of its
     \emph{gauged} parametrization $\delta U_1^\eta, \delta U_2^\eta, \delta
     U_3^\eta$ and $\delta \bS^\eta$, as
     in~\eqref{eq:tangenttuckerrepresentation} and~\eqref{eq:gauged_V_Tucker};
 \item the unknown $\textcolor{unknown}{\xi} \in {T_\bX \Mr}$ is to be
     determined 
 in terms of its \emph{gauged} parametrization $\textcolor{unknown}{\delta
     U_1}, \textcolor{unknown}{\delta U_2}, \textcolor{unknown}{\delta U_3}$
 and $\textcolor{unknown}{\delta \bS}$, again as
 in~\eqref{eq:tangenttuckerrepresentation} and~\eqref{eq:gauged_V_Tucker}.
\end{itemize}

To derive equations for $\textcolor{unknown}{\delta U_\mu}$ with $\mu = 1,2,3$ and $\textcolor{unknown}{\delta \bS}$ we exploit that
$T_\bX \Mr$ decomposes orthogonally into $\mcV_1 \oplus  \cdots \oplus
\mcV_{4}$; see~\eqref{eq:tangentdecomptucker}. This allows us to
split~\eqref{eq:GN_tucker}  into a system of four coupled equations by
projecting onto $\mcV_\mu$ for $\mu = 1,\ldots,4$. 

In particular, since $\textcolor{unknown}{\xi} \in T_\bX \Mr$ by assumption, we
can insert $\bZ := \mcL \Proj_{T_\bX \Mr} \textcolor{unknown}{\xi} = \mcL
\textcolor{unknown}{\xi}$ into~\eqref{eq:projectionTucker}. By exploiting the
structure of $\mcL$ (see~\eqref{eq:laplacelike}) and the orthogonality of the
gauged representation of tangent vectors (see~\eqref{eq:gauged_V_Tucker}), we
can simplify the expressions considerably and arrive at the equations
\begin{equation} \label{eq:coupledeqn}
 \begin{aligned}
    \dU_1^\eta &=\Proj_{U_1}^\perp \Big(L_1 U_1 \textcolor{unknown}{\delta S_{(1)}} + L_1 \textcolor{unknown}{\dU_1} S_{(1)} +
                \textcolor{unknown}{\dU_1} S_{(1)} \big[ I_{r_3} \otimes U_2^\T L_2 U_2 + U_3^\T L_3 U_3 \otimes I_{r_2} \big] 
                \Big) S_{(1)}^\dagger \\
    \dU_2^\eta &=\Proj_{U_1}^\perp \Big(L_2 U_2 \textcolor{unknown}{\delta S_{(2)}} + L_2 \textcolor{unknown}{\dU_2} S_{(2)} +
                \textcolor{unknown}{\dU_2} S_{(2)} \big[ I_{r_3} \otimes U_1^\T L_1 U_1 + U_3^\T L_3 U_3 \otimes I_{r_1} \big] 
                \Big) S_{(2)}^\dagger \\
    \dU_3^\eta &=\Proj_{U_1}^\perp \Big(L_3 U_3 \textcolor{unknown}{\delta S_{(3)}} + L_3 \textcolor{unknown}{\dU_3} S_{(3)} +
                \textcolor{unknown}{\dU_3} S_{(3)} \big[ I_{r_2} \otimes U_1^\T L_1 U_1 + U_2^\T L_2 U_2 \otimes I_{r_1} \big] 
                \Big) S_{(3)}^\dagger \\
    \delta \bS^\eta &= \big[ U_1^\T L_1 U_1 \textcolor{unknown}{\delta S_{(1)}} 
                                        + U_1^\T L_1 \textcolor{unknown}{\delta U_1} S_{(1)}\big]^{(1)} + \big[ U_2^\T L_2 U_2 \textcolor{unknown}{\delta S_{(2)}}  
                                        + U_2^\T L_2 \textcolor{unknown}{\delta U_2} S_{(2)}\big]^{(2)} \\
                                &\quad+ \big[ U_3^\T L_3 U_3  \textcolor{unknown}{\delta S_{(3)}}   
                                        + U_3^\T L_3 \textcolor{unknown}{\delta U_3} S_{(3)}\big]^{(3)}.
\end{aligned}
\end{equation}
Additionally, the gauge conditions need to be satisfied:
\begin{equation}\label{eq:dU_perp}
 U_1^\T \textcolor{unknown}{\delta U_1} = U_2^\T \textcolor{unknown}{\delta U_2} = U_3^\T \textcolor{unknown}{\delta U_3} = 0.
\end{equation}

In order to solve these equations, we will use the first three equations
of~\eqref{eq:coupledeqn}, together with~\eqref{eq:dU_perp}, to substitute
$\textcolor{unknown}{\delta U_\mu}$ in the last equation of~\eqref{eq:coupledeqn} and determine
a decoupled equation for $\textcolor{unknown}{\delta \bS}$. Rearranging the first equation
of~\eqref{eq:coupledeqn}, we obtain
\begin{align*}
\Proj_{U_1}^\perp \Big( L_1\textcolor{unknown}{\dU_1}  +
                \textcolor{unknown}{\dU_1} S_{(1)} \big[ I_{r_3} \otimes U_2^\T L_2 U_2 + U_3^\T L_3 U_3 \otimes I_{r_2} \big] 
                 S_{(1)}^\dagger \Big) =  \dU_1^\eta - \Proj_{U_1}^\perp L_1 U_1 \textcolor{unknown}{\delta S_{(1)}} S_{(1)}^\dagger.
\end{align*}
Vectorization and adhering to~\eqref{eq:dU_perp} yields the saddle point system
\begin{equation}\label{eq:saddle_Tucker}
    \begin{bmatrix}
        G
        & I_{r_1} \otimes U_1 \\
        I_{r_1} \otimes U_1^\T & {0}
    \end{bmatrix}
    \begin{bmatrix}
        \tenvec(\textcolor{unknown}{\dU_1})\\
        {y}_1 
    \end{bmatrix}
    =
    \begin{bmatrix}
       {b}_1\\
        {0}
    \end{bmatrix},
\end{equation}
where 
\begin{align*}
G &= I_{r_1} \otimes L_1 
        + (S_{(1)}^\dagger)^\T \big( I_{r_3} \otimes U_2^\T L_2 U_2 + U_3^\T L_3 U_3 \otimes I_{r_2} \big) S_{(1)}^\T \otimes I_{n_1},\\ 
		{b}_1 &= \tenvec(\dU^\eta_1) - \big( (S_{(1)}^\dagger)^\T \otimes
\Proj_{U_1}^\perp L_1 U_1 \big) \tenvec(\textcolor{unknown}{\delta S_{(1)}}),
\end{align*}
 and ${y}_1 \in \R^{r_1^2}$ is the dual variable. The positive definiteness of $L_1$ and the full rank
conditions on $U_1$ and $\bS$ imply that the above system is nonsingular; see, e.g.,
\cite{Benzi2005}. 
Using the Schur complement $G_S = -(I_{r_1} \otimes U_1)^\T G^{-1} (I_{r_1} \otimes U_1)$, 
we obtain the explicit expression
\begin{equation} \label{eq:expinverse}
 \tenvec(\textcolor{unknown}{\dU_1}) =  \Big( I_{n_1r_1} +  G^{-1}( I_{r_1} \otimes U_1) G_S^{-1} (I_{r_1} \otimes U_1^\T)\Big) G^{-1}  {b}_1 
 =  w_1 - F_1 \tenvec(\textcolor{unknown}{\delta S_{(1)}}), 
\end{equation} 
with
\begin{eqnarray*}
  w_1 &:=& \Big( I_{n_1 r_1} +  G^{-1} ( I_{r_1} \otimes U_1) G_S^{-1} (I_{r_1} \otimes U_1^\T)\Big) G^{-1}   \tenvec(\dU^\eta_1), \\
 F_1 &:=& \Big( I_{n_1 r_1} +  G^{-1}( I_{r_1} \otimes U_1) G_S^{-1} (I_{r_1} \otimes U_1^\T) \Big)G^{-1}   \Big( (S_{(1)}^\dagger)^\T \otimes \Proj_{U_1}^\perp L_1 U_1 \Big).
\end{eqnarray*}
Expressions analogous to~\eqref{eq:expinverse} can be derived for the other two factor matrices:
\begin{align*}
    \tenvec(\textcolor{unknown}{\dU_2}) &=  w_2 - F_2 \tenvec(\textcolor{unknown}{\delta S_{(2)}}), \\
    \tenvec(\textcolor{unknown}{\dU_3}) &=  w_3 - F_3 \tenvec(\textcolor{unknown}{\delta S_{(3)}}),
\end{align*}
with suitable analogs for $w_2, w_3, F_2$, and $F_3$. These expressions are now inserted into the last equation of~\eqref{eq:coupledeqn} for $\delta \bS^\eta$.
To this end, define permutation matrices $\Pi_{i \to j}$ that map the vectorization of the $i$th matricization to the vectorization of the $j$th matricization:
\[
       \Pi_{i \to j} \tenvec( \textcolor{unknown}{\delta \bS_{(i)}} ) = \tenvec( \textcolor{unknown}{\delta \bS_{(j)}} ),
\]
By definition, $\tenvec(\textcolor{unknown}{ \delta S_{(1)}} ) = \tenvec( \textcolor{unknown}{\delta \bS}
)$, and we finally obtain the following linear system
for $\tenvec(\textcolor{unknown}{\delta \bS})$:
\begin{equation} \label{eq:systemforS}
    F \tenvec( \textcolor{unknown}{\delta \bS})
                                 = \tenvec( \delta \bS^\eta ) - (S_{(1)}^\T
                                 \otimes U_1^\T L_1) w_1 - \Pi_{2\to 1}
                                 (S_{(2)}^\T \otimes U_2^\T L_2) w_2 -
                                 \Pi_{3\to 1} (S_{(3)}^\T \otimes U_3^\T L_3)
                                 w_3,
\end{equation}
with the $r_1r_2r_3\times r_1r_2 r_3$ matrix 
\begin{align*}
 F :=   I_{r_2r_3} & \otimes U_1^\T L_1 U_1 - (S_{(1)}^\T \otimes U_1^\T
 L_1)F_1  + \Pi_{2 \to 1}\Big[ I_{r_1r_3} \otimes U_2^\T L_2 U_2 - (S_{(2)}^\T
 \otimes U_2^\T L_2)F_2 \Big]\Pi_{1 \to 2} &\\
    &+\Pi_{3 \to 1}\Big[ I_{r_1r_2} \otimes U_3^\T L_3U_3 - (S_{(3)}^\T \otimes U_3^\T L_3)F_3 \Big]\Pi_{1 \to 3}.
\end{align*}
For small ranks, the linear system~\eqref{eq:systemforS} is solved by forming
the matrix $F$ explicitly and using a direct solver. Since this requires $O( r_1^3 r_2^3 r_3^3)$ operations, it is advisable to
use an iterative solver for larger ranks, in which the Kronecker product structure can be exploited when applying $F$; see also~\cite{Vandereycken2010}.
Once we have obtained $\textcolor{unknown}{\delta \bS}$, we can easily obtain
$\textcolor{unknown}{\delta U_1},\textcolor{unknown}{\delta
    U_2},\textcolor{unknown}{\delta U_3}$ using~\eqref{eq:expinverse}.

\begin{remark}\label{rem:Sylv_Tucker}
The application of $G^{-1}$ needed in~\eqref{eq:expinverse} as well as in the construction of
$G_S$ can be implemented
 efficiently by noting that $G$ is the matrix representation of the Sylvester operator
 $
  V \mapsto L_1 V + V \Gamma_1^\T,
 $
with the matrix
 \[
  \Gamma_1 := (S_{(1)}^\dagger)^\T \big( I_{r_3} \otimes U_2^\T L_2 U_2 + U_3^\T L_3 U_3 \otimes I_{r_2} \big) S_{(1)}^\T.
 \]
The $r_1 \times r_1$ matrix $\Gamma_1$ is non-symmetric but it can be diagonalized 
by first computing a QR decomposition $S_{(1)}^\T = Q_S R_S$ such that $Q_S^T Q_S = I_{r_1}$
and then computing the spectral decomposition of the symmetric matrix \[Q_S^\T \big( I_{r_3} \otimes U_2^\T L_2 U_2 + U_3^\T L_3 U_3 \otimes I_{r_2} \big) Q_S.\]
After diagonalization of $\Gamma_1$, the application of $G^{-1}$ requires the solution
of $r_1$ linear systems with the matrices $L_1 + \lambda I$, where $\lambda$ is an eigenvalue of $\Gamma_1$;
see also~\cite{Simoncini2013}. 
The Schur complement $G_S \in \R^{r_1^2 \times r_1^2}$ is constructed explicitly by applying $G^{-1}$ to the $r_1^2$ columns of  $I_{r_1} \otimes U_1$.

Analogous techniques apply to the computation of $w_2,F_2$, and $w_3,F_3$.
\end{remark}
Assuming, for example, that each $L_\mu$ is a tri-diagonal matrix, the solution of a linear system with the shifted
matrix $L_\mu + \lambda I$ can be performed in $O(n)$ operations. Therefore, using Remark~\ref{rem:Sylv_Tucker}, the
construction of the Schur complement $G_S$ requires $O(nr^3)$ operations. Hence,
the approximate Newton equation \eqref{eq:GN_tucker} can be solved in $O(dnr^3 + r^9)$ operations. This cost dominates the complexity of the Riemannian gradient calculation and the retraction step.

\subsection{The approximate Riemannian Hessian in the TT case}

When using the TT format, it seems to be much harder to solve
the approximate Newton equation~\eqref{eq:GN_manifold} directly and we therefore resort to
the preconditioned conjugate gradient (PCG) method for solving the linear system iteratively.
We use the following commonly used stopping criterion~\cite[Ch.~7.1]{Nocedal2006} for accepting the approximation $\tilde \xi$ produced by PCG:
\begin{equation*}
    \label{eq:newton_accuracy}
 \|  \Proj_{T_{\bX_k} \Mr} [\mcL  \tilde \xi - \nabla f(\bX_k)] \| \leq \min\Big( 0.5, \sqrt{\|  \Proj_{T_{\bX_k} \Mr} \nabla f(\bX_k) \|} \Big) \ \cdot \ \|  \Proj_{T_{\bX_k} \Mr} \nabla f(\bX_k) \|.
\end{equation*}
To derive an effective preconditioner for PCG, we first examine the approximate
Newton equation~\eqref{eq:GN_manifold}  more closely. For $d$-dimensional
tensors in the TT format, it takes the form
\begin{equation}\label{eq:GN_TT}
    \Proj_{T_\bX \Mr} \mcL \Proj_{T_\bX \Mr} \textcolor{unknown}{\xi} =  \eta,
\end{equation}
where
\begin{itemize}
    \item $\bX \in \Mr$ is parametrized by its cores $\bU_1,\bU_2,\ldots, \bU_d$ and is \emph{$d$-orthogonal};

    \item the right-hand side $\eta \in {T_\bX \Mr}$ is represented in terms of its \emph{gauged}
        parametrization $\delta \bU_1^\eta$, $\delta \bU_2^\eta$, $\ldots$, $\delta
        \bU_d^\eta$, as in~\eqref{eq:tangentTTrepresentation};
    \item the unknown $\textcolor{unknown}{\xi} \in {T_\bX \Mr}$ needs to be
        determined in terms of its \emph{gauged} parametrization $\textcolor{unknown}{\delta
        \bU_1}, \textcolor{unknown}{\delta \bU_2}, \ldots, \textcolor{unknown}{\delta \bU_d}$, again as in~\eqref{eq:tangentTTrepresentation}.
\end{itemize}
When PCG is applied to~\eqref{eq:GN_TT} with a preconditioner $\mcB\colon T_{\bX} \Mr \to T_{\bX} \Mr$, we need to evaluate
an expression of the form $\textcolor{unknown}{\xi} =  \mcB \eta$ for a given, arbitrary vector $\eta \in {T_\bX \Mr}$.
Again, $\textcolor{unknown}{\xi}$ and $\eta$ are represented using the gauged parametrization above.

We will present two block Jacobi preconditioners for~\eqref{eq:GN_TT}; both are
variants of parallel subspace correction (PSC) methods~\cite{Xu:1992}. They mainly differ in the way the tangent space $T_\bX
\Mr$ is split into subspaces. 

\subsubsection{A block diagonal Jacobi preconditioner}
\label{subsubsec:blockjacobi}

The most immediate choice for splitting $T_\bX \Mr$ is to simply take the
direct sum~\eqref{eq:tangentdecomp}. The  PSC method is then defined in terms
of the local operators 
\[
 \mcL_\mu \colon \mcV_\mu \to \mcV_\mu, \qquad \mcL_\mu = \left.\Proj^\mu \mcL \Proj^\mu\right|_{\mcV_\mu}, \qquad \mu=1,\ldots,d,
\]
where $\Proj^\mu$ is the orthogonal projector onto $\mcV_\mu$; see
\S\ref{sec:intro_TT}. The operators $\mcL_\mu$ are symmetric and positive
definite, and hence invertible, on $\mcV_\mu$. This allows us to express the
resulting preconditioner as \cite[\S3.2]{Xu2001}
\[
 \mcB = \sum_{\mu=1}^d  \mcL_\mu^{-1} \Proj^\mu = \sum_{\mu=1}^d \left(\left.\Proj^\mu \mcL \Proj^\mu\right|_{\mcV_\mu}\right)^{-1} \Proj^\mu .
\]
The action of the preconditioner $\textcolor{unknown}{\xi} = \mcB \eta$ can thus be computed as $\textcolor{unknown}{\xi} = \sum_{\mu=1}^d \textcolor{unknown}{\xi_\mu}$ with
\[
 \textcolor{unknown}{\xi_\mu} = \left(\left.\Proj^\mu \mcL \Proj^\mu\right|_{\mcV_\mu}\right)^{-1} \Proj^\mu \eta, \qquad \mu = 1,\ldots, d.
\]
\begin{paragraph}{Local problems.}
 The local equations determining $\textcolor{unknown}{\xi_\mu}$,
\begin{equation}\label{eq:local_solve_nonverlap}
 \Proj^\mu \mcL \Proj^\mu \textcolor{unknown}{\xi_\mu} = \Proj^\mu \eta, \qquad \textcolor{unknown}{\xi_\mu} \in \mcV_\mu, \qquad \mu = 1,\ldots, d,
\end{equation}
can be solved for all $\textcolor{unknown}{\xi_\mu} \in \mcV_\mu$ in parallel.
By~\eqref{eq:V_TT}, we have $\textcolor{unknown}{\xi_\mu} = \bX_{\neq \mu}
\tenvec(\textcolor{unknown}{\delta \bU_\mu})$ for some gauged
$\textcolor{unknown}{\delta \bU_\mu}$. Since $\Proj^\mu \eta$ satisfies an
expansion analogous to~\eqref{eq:tangentTTrepresentation}, straightforward
properties of the projectors $\Proj^\mu$ allow us to
write~\eqref{eq:local_solve_nonverlap} as
\begin{align*}
    \quad \Proj^\mu \mcL \bX_{\neq \mu} \tenvec(\textcolor{unknown}{\delta \bU_\mu}) &= \bX_{\neq \mu} \tenvec(\delta \bU^\eta_\mu), \qquad \mu = 1,\ldots, d,
\end{align*}
under the additional constraint $(\textcolor{unknown}{\delta \bU^\unl_\mu})^\T \bU_\mu^\unl = 0$ when  $\mu \neq d$. Now expressing
the result of $\Proj^\mu$ applied to $\mcL \bX_{\neq \mu} \tenvec(\textcolor{unknown}{\delta \bU_\mu})$ as in~\eqref{eq:P_Z_TT_in_components} and using~\eqref{eq:deltaU} leads to 
\begin{equation}
    \label{eq:muthcomptangent}
    (I_{n_\mu r_{\mu-1}} - \Proj^\unl_\mu)\big(I_{n_\mu} \otimes \bX_{\leq
        \mu-1}^\T\big) \Big[\mcL \bX_{\neq \mu}
    \tenvec(\textcolor{unknown}{\delta \bU_\mu})\Big]^{<\mu>} \bX_{\geq
        \mu+1}\big( \bX_{\geq \mu+1}^\T \bX_{\geq \mu + 1}\big)^{-1} =
	(\delta\bU_\mu^\eta)^\unl
\end{equation}
for $\mu\neq d$, while~\eqref{eq:deltaU_d} for $\mu=d$ leads to  the equation
\begin{equation}
    \label{eq:muthcomptangent_d}
    \big(I_{n_d} \otimes \bX_{\leq d-1}^\T\big) \Big[\mcL \bX_{\neq d} \tenvec(\textcolor{unknown}{\delta \bU_d})\Big]^{<d>} 
	= (\delta\bU_d^\eta)^\unl.
\end{equation}

Using~\eqref{eq:Xneqmu}, the application of the Laplace-like operator $\mcL$ to $\bX_{\neq \mu}$ can be decomposed into three parts,
\begin{equation}\label{eq:L_neq_mu_onesided}
    \mcL \bX_{\neq \mu} = \widetilde{\mcL}_{\geq \mu+1} \otimes I_{n_\mu} \otimes \bX_{\leq\mu-1} +
    \bX_{\geq\mu+1} \otimes L_\mu \otimes \bX_{\leq\mu-1}+
    \bX_{\geq\mu+1} \otimes I_{n_\mu} \otimes \widetilde{\mcL}_{\leq \mu-1} 
\end{equation}
with the reduced leading and trailing terms
\[
    \widetilde{\mcL}_{\leq \mu-1} = \left(\sum_{\nu = 1}^{\mu-1} I_{n_{\mu-1}}
    \otimes \dots \otimes L_{\nu} \otimes \dots \otimes I_{n_1} \right)\bX_{\le
    \mu-1}^{},
\]
\[
    \widetilde{\mcL}_{\geq \mu+1} = \left(\sum_{\nu = \mu+1}^d I_{n_d} \otimes
    \dots \otimes L_\nu \otimes \dots I_{n_{\mu+1}} \right)\bX_{\ge \mu+1}^{}.
\]
Some manipulation\footnote{This is shown by applying the relation $\bX^{<\mu>}
    = (I_{n_\mu} \otimes \bX_{\leq \mu-1}) \bU^\unl_\mu \bX_{\geq \mu+1}^\T$,
    which holds for any TT tensor~\cite[eq.~(2.4)]{Lubich2014}, to $\mcL
    \bX_{\neq \mu}\tenvec(\textcolor{unknown}{\delta \bU_\mu})$.} establishes
the identity
\[
    \big[\mcL \bX_{\neq \mu} \tenvec(\textcolor{unknown}{\delta \bU_\mu})\big]^{<\mu>} = 
     \left( I_{n_\mu}\otimes \bX_{\leq\mu-1}\right)\textcolor{unknown}{\delta \bU_\mu^\unl} \widetilde{\mcL}_{\geq \mu+1}^\T +
            \left(L_\mu \otimes \bX_{\leq\mu-1}+
                I_{n_\mu} \otimes \widetilde{\mcL}_{\leq \mu-1}\right) \textcolor{unknown}{\delta \bU_\mu^\unl} \bX_{\geq\mu+1}^\T.
\]
Inserting this expression into \eqref{eq:muthcomptangent} yields for $\mu \neq d$
\begin{multline*}    
   (I_{n_\mu r_{\mu-1}} - \Proj^\unl_\mu) \big[      \textcolor{unknown}{\delta \bU_\mu^\unl} \widetilde{\mcL}_{\geq \mu+1}^\T
         \bX_{\geq \mu+1}\big( \bX_{\geq \mu+1}^\T \bX_{\geq \mu + 1}\big)^{-1}
        \\ +(L_\mu \otimes I_{r_{\mu-1}}+
                I_{n_\mu} \otimes \bX_{\leq \mu-1}^\T\widetilde{\mcL}_{\leq \mu-1}) \; \textcolor{unknown}{\delta \bU_\mu^\unl} \big]
				 = (\delta\bU_\mu^\eta)^\unl.
\end{multline*}
After defining the (symmetric positive definite) matrices $\mcL_{\leq \mu-1} = \bX_{\leq \mu-1}^\T\widetilde{\mcL}_{\leq
\mu-1}$ and $\mcL_{\geq \mu+1} = \bX_{\geq
\mu+1}^\T \widetilde{\mcL}_{\geq \mu+1}$, we finally obtain
\begin{equation}
    \label{eq:jacobieq}
         (I_{n_\mu r_{\mu-1}} - \Proj^\unl_\mu) \big[ \textcolor{unknown}{\delta \bU_\mu^\unl} \mcL_{\geq \mu+1}
         \big( \bX_{\geq \mu+1}^\T \bX_{\geq \mu + 1}\big)^{-1}
         +(L_\mu \otimes I_{r_{\mu-1}}+I_{n_\mu} \otimes \mcL_{\leq \mu-1})  \textcolor{unknown}{\delta \bU_\mu^\unl} \big]
				 = (\delta\bU_\mu^\eta)^\unl,
\end{equation}
with the gauge condition $(\textcolor{unknown}{\delta \bU^\unl_\mu})^\T \bU_\mu^\unl = 0$. For $\mu = d$, there is no gauge condition and~\eqref{eq:muthcomptangent_d} becomes
\begin{align}\label{eq:jacobieq_d}
         \textcolor{unknown}{\delta \bU_d^\unl} 
         +\left(L_d \otimes I_{r_d}+
                I_{n_d} \otimes \mcL_{\leq d-1}\right) \textcolor{unknown}{\delta \bU_d^\unl}
         =     (\delta\bU_d^\eta)^\unl.				
\end{align}

\end{paragraph}

\begin{paragraph}{Efficient solution of local problems.}
The derivations above have led us to the linear systems~\eqref{eq:jacobieq}
and~\eqref{eq:jacobieq_d} for determining the local component
$\textcolor{unknown}{\xi_\mu}$. While~\eqref{eq:jacobieq_d} is a Sylvester
equation and can be solved with standard techniques, more work is needed to
address~\eqref{eq:jacobieq} efficiently.
Since $\mcL_{\geq \mu+1}$ and $ \bX_{\geq \mu+1}^\T \bX_{\geq \mu + 1}$ are
symmetric positive definite, they admit a generalized eigenvalue decomposition:
There is an invertible matrix $Q$ such that $\mcL_{\geq \mu+1} Q =  (\bX_{\geq
    \mu+1}^\T \bX_{\geq \mu + 1}) Q \Lambda$ with $\Lambda$ diagonal and $Q^\T
( \bX_{\geq \mu+1}^\T \bX_{\geq \mu + 1}) Q = I_{r_\mu}$. This
transforms~\eqref{eq:jacobieq} into
\begin{equation*}
         (I_{n_\mu r_{\mu-1}} - \Proj^\unl_\mu) \big[ \textcolor{unknown}{\delta \bU_\mu^\unl} Q^\T \Lambda
         +\left(L_\mu \otimes I_{r_\mu}+
                I_{n_\mu} \otimes \mcL_{\leq \mu-1}\right) \textcolor{unknown}{\delta \bU_\mu^\unl} Q^\T \big]
				=     (\delta\bU_\mu^\eta)^\unl  Q^\T.
\end{equation*}
Setting 
 $\textcolor{unknown}{\widetilde{\delta\bU}{}_\mu^\unl} =\textcolor{unknown}{\delta \bU_\mu^\unl} Q^\T $ and 
 $(\widetilde{\delta\bU}{}^\eta_\mu)^\unl =(\delta\bU_\mu^\eta)^\unl Q^\T $, we can formulate these equations column-wise:
\begin{equation} \label{eq:columnwisesomething}
    \begin{aligned}
 (I_{n_\mu r_{\mu-1}} - \Proj^\unl_\mu) \left[\lambda_i  I_{r_\mu n_\mu} 
         + L_\mu \otimes I_{r_\mu}+
        I_{n_\mu} \otimes \mcL_{\leq \mu-1} \right] \textcolor{unknown}{\widetilde{\delta\bU}{}_\mu^\unl(:,i)} 
  =  (\widetilde{\delta\bU}{}^\eta_\mu)^\unl(:,i) ,
    \end{aligned}
\end{equation}
where $\lambda_i = \Lambda(i,i) > 0$. 
Because $Q$ is invertible, the gauge-conditions on $\textcolor{unknown}{\delta \bU_\mu^\unl}$ are equivalent to $(\textcolor{unknown}{\widetilde{\delta\bU}{}_\mu^\unl})^\T \bU_\mu^\unl = 0$.
Combined with~\eqref{eq:columnwisesomething}, we obtain -- similar to~\eqref{eq:saddle_Tucker} -- the saddle point systems
\begin{equation}
    \label{eq:jacobi_saddle}
    \begin{bmatrix}
        G_{\mu,i}
        & \bU_\mu^\unl \\
        (\bU_\mu^\unl)^\T & {0}
    \end{bmatrix}
    \begin{bmatrix}
        \textcolor{unknown}{\widetilde{\delta\bU}{}_\mu^\unl(:,i)} \\
        y
    \end{bmatrix}
    =
    \begin{bmatrix}
        (\widetilde{\delta\bU}{}^\eta_\mu)^\unl(:,i) \\ 0
    \end{bmatrix}
\end{equation}
with the symmetric positive definite matrix
 \begin{equation} \label{eq:defgmui}
 G_{\mu,i} = \lambda_i  I_{n_\mu} \otimes I_{r_\mu} + L_\mu \otimes I_{r_\mu}+ I_{n_\mu} \otimes \mcL_{\leq \mu-1}
\end{equation}
 and the dual variable $y \in \R^{r_\mu}$. The system~\eqref{eq:jacobi_saddle} is solved for each column of $\textcolor{unknown}{\widetilde{\delta\bU}{}_\mu^\unl}$:
\[
 \textcolor{unknown}{\widetilde{\delta\bU}{}_\mu^\unl(:,i)} =  \Big( I_{n_\mu r_\mu} +  G_{\mu,i}^{-1} \; \bU_\mu^\unl \; G_S^{-1} \; (\bU_\mu^\unl)^\T \Big) G_{\mu,i}^{-1}  \; (\widetilde{\delta\bU}{}^\eta_\mu)^\unl(:,i) ,
\]
using the Schur complement $G_S := - (\bU_\mu^\unl)^\T G_{\mu,i}^{-1} \bU_\mu^\unl$. Transforming back eventually yields   $\textcolor{unknown}{\delta\bU_\mu^\unl} = \textcolor{unknown}{\widetilde{\delta\bU}{}_\mu^\unl} Q^{-\T}$.

\begin{remark}\label{rem:Sylv_TT}
Analogous to Remark~\ref{rem:Sylv_Tucker}, the application of $G_{\mu,i}^{-1}$ benefits from the fact
 that the matrix $G_{\mu,i}$ defined in~\eqref{eq:defgmui} represents the Sylvester operator
 $$
  V \mapsto (L_\mu + \lambda_i I_{n_\mu}) V + V \mcL_{\leq \mu-1}.
 $$
 After diagonalization of $\mcL_{\leq \mu-1}$, the application of $G_{\mu,i}^{-1}$ requires the solution
of $r_\mu$ linear systems with the matrices $L_\mu + (\lambda_i + \beta) I_{n_\mu}$, where $\beta$ is an eigenvalue of $\mcL_{\leq \mu-1}$. 
The Schur complements $G_S \in \R^{r_\mu \times r_\mu}$ are constructed explicitly by applying $G_{\mu,i}^{-1}$ to the $r_\mu$ columns of  $\bU_\mu^\unl$. 
\end{remark}
Assuming again that solving with the shifted matrices $L_\mu + (\lambda_i + \beta)
I_{n_\mu}$ can be performed in $O(n_\mu)$ operations, the construction of the
Schur complement $G_S$ needs $O(n_\mu r_\mu^2)$ operations. Repeating this for
all $r_\mu$ columns of $\textcolor{unknown}{\widetilde{\delta\bU}{}_\mu^\unl}$
and all cores $\mu = 1,\ldots,d-1$ yields a total computational complexity of $O(dnr^3)$ for applying the block-Jacobi preconditioner.
\end{paragraph}

\subsubsection{An overlapping block-Jacobi preconditioner} \label{sec:overlapping}
The block diagonal preconditioner discussed above is computationally expensive
due to the need for solving the saddle point systems~\eqref{eq:jacobi_saddle}.
To avoid them, we will construct a PSC preconditioner for the subspaces
\[
    \widehat{\mcV}_\mu := \left\{ \,\bX_{\neq \mu} \tenvec(\delta \bU_\mu) \colon \delta \bU_\mu \in \R^{r_{\mu-1} \times n_\mu \times r_\mu} \right\}  = \myspan \bX_{\neq \mu}, \qquad \mu = 1, \ldots, d.
\]
Observe that $\mcV_\mu \subsetneq \widehat{\mcV}_\mu$ for $\mu \neq d$. Hence,
the decomposition $T_\bX \Mr = \cup_{\mu=1}^d \widehat{\mcV}_\mu $ is no longer
a direct sum as in~\eqref{eq:tangentdecomp}. The advantage of
$\widehat{\mcV}_\mu$ over $\mcV_\mu$, however, is that the orthogonal projector
$\widehat{\Proj}{}^\mu$ onto $\widehat{\mcV}_\mu$ is considerably easier. In
particular, since $\bX$ is $d$-orthogonal, we obtain
\begin{equation}\label{eq:Pmu_hat}
    \widehat{\Proj}{}^\mu =  \bX_{\neq \mu} (\bX_{\neq \mu}^\T \bX_{\neq \mu})^{-1} \bX_{\neq \mu}^\T = 
	\bX_{\neq \mu} \left[ (\bX_{\ge {\mu+1}}^\T \bX_{\ge {\mu+1}})^{-1} \otimes I_{n_\mu} \otimes I_{r_{\mu-1}} \right] \bX_{\neq \mu}^\T .
\end{equation}

The PSC preconditioner corresponding to the subspaces $\widehat{\mcV}_\mu$ is given by
\[
 \widehat\mcB = \sum_{\mu=1}^d \left(\left.\widehat{\Proj}{}^\mu \mcL \widehat{\Proj}{}^\mu\right|_{\widehat{\mcV}_\mu}\right)^{-1} \widehat{\Proj}{}^\mu.
\]
The action of the preconditioner $\textcolor{unknown}{\xi} = \widehat \mcB \eta$ can thus be computed as $\textcolor{unknown}{\xi} = \sum_{\mu=1}^d \textcolor{unknown}{\xi_\mu}$ with
\begin{equation}\label{eq:local_solve_oververlap}
 \widehat{\Proj}{}^\mu \mcL \widehat{\Proj}{}^\mu \textcolor{unknown}{\xi_\mu} = \widehat{\Proj}{}^\mu \eta, \qquad \textcolor{unknown}{\xi_\mu} \in \widehat{\mcV}_\mu, \qquad \mu = 1,\ldots, d.
\end{equation}

\begin{paragraph}{Local problems.}
To solve the local equations~\eqref{eq:local_solve_oververlap}, we proceed as
in the previous section, but the resulting equations will be considerably
simpler.  Let $\widehat{\Proj}{}^\mu \eta = \bX_{\neq \mu}
\tenvec(\widehat{\delta\bU}{}_\mu^\eta)$ for some
$\widehat{\delta\bU}{}_\mu^\eta$, which will generally differ from the gauged
$\delta\bU_\mu^\eta$ parametrization of $\eta$. Writing
$\textcolor{unknown}{\xi_\mu} = \bX_{\neq \mu}
\tenvec(\textcolor{unknown}{\delta \bU_\mu})$, we obtain the linear systems
\[
 \widehat{\Proj}{}^\mu \mcL \bX_{\neq \mu} \tenvec(\textcolor{unknown}{\delta \bU_\mu}) = \bX_{\neq \mu} \tenvec(\widehat{\delta\bU}{}_\mu^\eta)
\]
for $\mu = 1,\ldots, d$.
Plugging in~\eqref{eq:Pmu_hat} gives
\begin{equation}
    \label{eq:jacobieq_ten}
\left[ (\bX_{\ge {\mu+1}}^\T \bX_{\ge {\mu+1}})^{-1} \otimes I_{n_\mu} \otimes
    I_{r_{\mu-1}} \right] \bX_{\neq \mu}^\T \mcL \bX_{\neq \mu}
\tenvec(\textcolor{unknown}{\delta \bU_\mu}) =
\tenvec(\widehat{\delta\bU}{}_\mu^\eta).
\end{equation}
Analogous to~\eqref{eq:L_neq_mu_onesided}, we can write 
\[
    \bX_{\neq \mu}^\T \mcL \bX_{\neq \mu} = {\mcL}_{\geq \mu+1} \otimes I_{n_\mu} \otimes I_{r_{\mu-1}} +
    \bX_{\geq\mu+1}^\T \bX_{\geq\mu+1} \otimes L_\mu \otimes I_{r_{\mu-1}}+
    \bX_{\geq\mu+1}^\T \bX_{\geq\mu+1} \otimes I_{n_\mu} \otimes {\mcL}_{\leq \mu-1} 
\]
with the left and right parts
\[
    {\mcL}_{\leq \mu-1} = \bX_{\le
    \mu-1}^{\T} \left(\sum_{\nu = 1}^{\mu-1} I_{n_{\mu-1}}
    \otimes \dots \otimes L_{\nu} \otimes \dots \otimes I_{n_1} \right)\bX_{\le
    \mu-1}^{},
\]
\[
    {\mcL}_{\geq \mu+1} = \bX_{\ge \mu+1}^{\T} \left(\sum_{\nu = \mu+1}^d I_{n_d} \otimes
    \dots \otimes L_\nu \otimes \dots I_{n_{\mu+1}} \right)\bX_{\ge \mu+1}^{}.
\]
Again, it is not hard to show that
\[
     \big(\bX_{\neq \mu}^\T\mcL \bX_{\neq \mu}\tenvec(\textcolor{unknown}{\delta \bU_\mu})\big)^{<\mu>} = 
\textcolor{unknown}{\delta \bU_\mu^\unl} {\mcL}_{\geq \mu+1} +
            \left(L_\mu \otimes  I_{r_{\mu-1}}+
                I_{n_\mu} \otimes {\mcL}_{\leq \mu-1}\right) \textcolor{unknown}{\delta \bU_\mu^\unl} \bX_{\geq\mu+1}^\T  \bX_{\geq\mu+1}.
\]
Hence, \eqref{eq:jacobieq_ten} can be written as
\begin{equation}
    \label{eq:jacobieq_unl}
\textcolor{unknown}{\delta \bU_\mu^\unl} \mcL_{\geq \mu+1}
         \big( \bX_{\geq \mu+1}^\T \bX_{\geq \mu + 1}\big)^{-1}
         +(L_\mu \otimes I_{r_{\mu-1}}+I_{n_\mu} \otimes \mcL_{\leq \mu-1})  \textcolor{unknown}{\delta \bU_\mu^\unl} 
				 = (\widehat{\delta\bU}{}_\mu^\eta)^\unl.
\end{equation}
\end{paragraph}
\begin{paragraph}{Efficient solution of local problems.}
The above equations can be directly solved as follows: Using the generalized
eigendecomposition of $\mcL_{\geq \mu+1} Q =  (\bX_{\geq \mu+1}^\T \bX_{\geq
    \mu + 1}) Q \Lambda$, we can write~\eqref{eq:jacobieq_unl} column-wise as
\[
 G_{\mu,i} \; \textcolor{unknown}{\widetilde{{\delta\bU}}{}_\mu^\unl(:,i)}  = (\widetilde{\widehat{\delta\bU}}{}^\eta_\mu)^\unl(:,i)
\]
with the system matrix
 \[
 G_{\mu,i} = \lambda_i  I_{n_\mu} \otimes I_{r_\mu} + L_\mu \otimes I_{r_\mu}+ I_{n_\mu} \otimes \mcL_{\leq \mu-1}, \qquad \lambda_i = \Lambda(i,i),
\]
and the transformed variables
$\textcolor{unknown}{\widetilde{{\delta\bU}}{}_\mu^\unl}
:=\textcolor{unknown}{\delta \bU_\mu^\unl} Q^\T $ and
$(\widetilde{\widehat{\delta\bU}}{}^\eta_\mu)^\unl :=
(\widehat{\delta\bU}{}_\mu^\eta)^\unl Q^\T $.
  Solving with $G_{\mu,i}$ can again be achieved by efficient solvers for
  Sylvester equations, see Remark~\ref{rem:Sylv_TT}. After forming
  $\textcolor{unknown}{\delta \bU_\mu^\unl} =
  \textcolor{unknown}{\widetilde{{\delta\bU}}{}_\mu^\unl} Q^{-\T}$ for all
  $\mu$, we have obtained the solution as an ungauged parametrization:
\[
    \textcolor{unknown}{\xi} = \widehat \mcB \eta = \sum_{\mu=1}^d \bX_{\neq \mu} \tenvec(\textcolor{unknown}{\delta \bU_\mu}) .
 \]
 To obtain the gauged parametrization of $\textcolor{unknown}{\xi}$
 satisfying~\eqref{eq:tangentTTrepresentation}, we can simply
 apply~\eqref{eq:deltaU} to compute $
 \Proj_{T_{\bX}\Mr}(\textcolor{unknown}{\xi})$ and exploit that
 $\textcolor{unknown}{\xi}$ is a TT tensor (with doubled TT ranks compared to
 $\bX$).

Assuming again that solving with $L_\mu$ can be performed in $O(n_\mu)$
operations, we end up with a total computational complexity of $O(dnr^3)$ for
applying the overlapping block-Jacobi preconditioner. Although this is the same
asymptotic complexity as the non-overlapping scheme from \S
\ref{subsubsec:blockjacobi}, the constant and computational time can be
expected to be significantly lower thanks to not having to solve saddle point
systems in each step.

\begin{remark}
    By $\mu$-orthogonalizing $\bX$ and transforming $\textcolor{unknown}{\delta
        \bU_\mu}$, as described in~\cite{Steinlechner2015}, the Gram matrix
    $\bX_{\geq \mu+1}^\T \bX_{\geq \mu + 1}$ in~\eqref{eq:jacobieq}
    and~\eqref{eq:jacobieq_unl} becomes the identity matrix. This leads to a
    more stable calculation of the corresponding unknown
    $\textcolor{unknown}{\delta \bU_\mu}$, see also Remark~\ref{rem:inv_gram}.
    We make use of this transformation in our implementations.
\end{remark}	
\end{paragraph}

\subsubsection{Connection to ALS}
The overlapping block-Jacobi preconditioner $\widehat \mcB$ introduced above is
closely related to ALS applied to~\eqref{eq:calAxf}. There are, however,
crucial differences explaining why $\widehat \mcB$ is significantly cheaper per
iteration than ALS.

Using $\tenvec(\bX) = \bX_{\neq \mu} \tenvec(\bU_\mu)$, one micro-step of ALS fixes $\bX_{\neq \mu}$ and replaces $\bU_\mu$ by the minimizer of (see, e.g., \cite[Alg.~1]{Holtz2012})
\[
 \min_{\bU_\mu} \frac{1}{2} \langle \bX_{\neq \mu} \tenvec(\bU_\mu), \mcA \bX_{\neq \mu} \tenvec(\bU_\mu) \rangle - \langle \bX_{\neq \mu} \tenvec(\bU_\mu), \tenvec(\bF) \rangle.
\]
After $\bU_\mu$ has been updated, ALS proceeds to the next core until all cores
have eventually been updated in a particular order, for example, $\bU_1, \bU_2,
\ldots, \bU_d$.  The solution to the above minimization problem is obtained
from solving the ALS subproblem
\[
 \bX_{\neq \mu}^\T \mcA \bX_{\neq \mu} \tenvec(\bU_\mu)  = \bX_{\neq \mu}^\T \tenvec(\bF) .
\]

 It is well-known that ALS can be seen as a block version of non-linear
 Gauss--Seidel. The subproblem typically needs to be computed iteratively since
 the system matrix $\bX_{\neq \mu}^\T \mcA \bX_{\neq \mu} \bU_\mu$ is often
 unmanageably large.

When $\bX$ is $\mu$-orthogonal, $\bX_{\ge {\mu+1}}^\T \bX_{\ge {\mu+1}} = I_{r_\mu}$ and the ALS subproblem has the same form as the subproblem~\eqref{eq:jacobieq_ten} in the overlapping block-Jacobi preconditioner $\widehat \mcB$. However, there are crucial differences:
\begin{itemize}
    \item ALS directly optimizes for the cores and as such uses $\mcA$ in the
        optimization problem. The approximate Newton method, on the other hand,
        updates (all) the cores using a search direction obtained from
        minimizing the quadratic model~\eqref{eq:GN_manifold}. It can therefore
        use any positive definite approximation of $\mcA$ to construct this
        model, which we choose as  $\mcL$. Since~\eqref{eq:jacobieq_ten} is the
        preconditioner for this quadratic model, it uses $\mcL$ as well.
	\item ALS updates each core immediately and it is a block version of non-linear Gauss--Seidel for~\eqref{eq:calAxf}, whereas $\widehat \mcB$ updates all the cores simultaneously resembling a block version of linear Jacobi. 	
    \item Even in the  large-scale setting of $n_\mu \gg 10^3$, the
        subproblems~\eqref{eq:jacobieq_ten}  can be solved efficiently in
        closed form as long as $L_\mu+\lambda I_{n_\mu}$ allows for efficient
        system solves, e.g., for tridiagonal $L_\mu$. This is not possible in
        ALS where the subproblems have to be formulated with $\mcA$ and
        typically need to be solved iteratively using PCG.
\end{itemize}
\begin{remark}
	Instead of PSC, we experimented with a symmetrized version of a successive subspace correction (SSC) preconditioner, also known as a back and forth ALS sweep. However, the higher computational cost per iteration of SSC was not offset by a possibly improved convergence behavior.
\end{remark}

\section{Numerical experiments}
\label{sec:experiments}
In this section, we compare the performance of the different preconditioned optimization techniques
discussed in this paper for two representative test cases.

We have implemented all algorithms in \textsc{Matlab}. For the TT format, we have made use of the TTeMPS toolbox, see \texttt{http://anchp.epfl.ch/TTeMPS}.
All numerical experiments and timings are performed on a 12 core Intel Xeon
X5675, 3.07 Ghz, 192 GiB RAM using \textsc{Matlab} 2014a, running on
Linux kernel 3.2.0-0.

To simplify the discussion, we assume throughout this section that the tensor
size and ranks are equal along all modes and therefore state them as scalar
values: $n = \max_\mu n_\mu$ and $r = \max_\mu r_\mu$.

\subsection{Test case 1: Newton potential}
\label{subsec:newtonpotential}
As a standard example leading to a linear system of the form \eqref{eq:laplace+potential}, we consider the partial differential equation 
\begin{alignat*}{2}
    -\Delta u(x) + V(x) &= f(x),\quad  &&x \in \Omega = (-10,10)^d,\\
    u(x) &= 0                              &&x \in \partial \Omega. 
\end{alignat*}
with the Laplace operator $\Delta$, the \emph{Newton potential} $V(x) =
\|x\|^{-1}$, and the source function $f\colon\R^d \to \R$. Equations of this type are used
to describe the energy of a charged particle in an electrostatic potential. 

We discretize the domain $\Omega$ by a uniform tensor grid with $n^d$ grid points and corresponding mesh width $h$.
Then, by finite difference approximation on this tensor grid, we obtain a tensor equation
of the type \eqref{eq:calAxf}, where the linear operator $\mcA$ is the sum of 
the $d$-dimensional Laplace operator as in \eqref{eq:laplacelike} with $L_\mu = \frac{1}{h^2}\operatorname{tridiag}(-1, 2, -1)\in \R^{n \times n}$, and the discretized Newton potential $\bV$.
To create a low-rank representation of the Newton potential, $V(x)$ is approximated by a
rank 10 tensor $\bV$ using exponential sums~\cite{Hackbusch2010}. The application of $\mcA$ to a 
tensor $\bX$ is given by
\[
    \mcA \bX = \mcL \bX + \bV \circ \bX,
\]
where $\circ$ denotes the Hadamard (element-wise) product. The application of this
operator increases the ranks significantly: If $\bX$ has rank $r$
then $\mcA \bX$ has rank $(2 + 10) r = 12 r$.

\subsection{Test case 2: Anisotropic Diffusion Equation}
\label{subsec:diffusion}
As a second example, we consider the anisotropic diffusion equation
\begin{alignat*}{2}
    -\operatorname{div} (D \nabla u(x)) &= f(x),\quad  &&x \in \Omega = (-10,10)\,^d,\\
    u(x) &= 0                              &&x \in \partial \Omega,
\end{alignat*}
with a tridiagonal diffusion matrix
$D = \operatorname{tridiag}(\alpha, 1, \alpha)\in \R^{d \times d}$.
The discretization on a uniform tensor grid with $n^d$ grid points and mesh width
$h$ yields a linear equation with system matrix $A = L + V$ consisting of the potential term 
\[
    V = I_{n} \otimes \cdots \otimes I_n\otimes B_2 \otimes 2 \alpha B_1
    \,+\, I_{n} \otimes \cdots \otimes I_n \otimes B_3 \otimes 2 \alpha B_2 \otimes I_n 
    \,+\, B_d \otimes 2\alpha B_{d-1} \otimes I_n \otimes \cdots \otimes I_n,
\]
and the Laplace part $L$ defined as in the previous example. The matrix $B_\mu =
\frac{1}{2h}\operatorname{tridiag}(-1, 0, 1)\in \R^{n \times n}$ represents the
one-dimensional central finite difference matrix for the first derivative.

\begin{figure}[h]
    \centering
    \begin{tikzpicture}[scale=0.5]
        \definecolor{blue}{RGB}{0,106,214}
        \draw[draw=blue!50] (0,0,0) -- (3,0,0) -- (3,3,0) -- (0,3,0) -- cycle;
        \draw[draw=blue!50] (0,0,0) -- (3,0,0) -- (3,0,3) -- (0,0,3) -- cycle;
        \begin{scope}[shift={(1,0,0)}]
            \fill[blue!10] (0,1,0) -- (1,1,0) -- (1,1,1) -- (0,1,1) -- cycle;
            \shade[bottom color=blue!60, top color=blue!30] (1,0,0) -- (1,1,0) -- (1,1,1) -- (1,0,1) -- cycle;
            \shade[bottom color=blue!60, top color=blue!30] (0,0,1) -- (1,0,1) -- (1,1,1) -- (0,1,1) -- cycle;
            \draw[blue!70, dashed] (1,0,0) -- (0,0,0) -- (0,1,0);
            \draw[blue!70, dashed] (0,0,0) -- (0,0,1);
            \draw[blue] (0,1,1) -- (0,0,1) -- (1,0,1) -- (1,1,1);
            \draw[blue] (1,0,1) -- (1,0,0) -- (1,1,0);
            \path (0.35,0.35,0.5) node[scale=1,black] {$\bS$};
        \end{scope}
        \begin{scope}[shift={(0,1,1)}]
            \fill[blue!10] (0,1,0) -- (1,1,0) -- (1,1,1) -- (0,1,1) -- cycle;
            \shade[bottom color=blue!60, top color=blue!30] (1,0,0) -- (1,1,0) -- (1,1,1) -- (1,0,1) -- cycle;
            \shade[bottom color=blue!60, top color=blue!30] (0,0,1) -- (1,0,1) -- (1,1,1) -- (0,1,1) -- cycle;
            \draw[blue!70, dashed] (1,0,0) -- (0,0,0) -- (0,1,0);
            \draw[blue!70, dashed] (0,0,0) -- (0,0,1);
            \draw[blue] (0,1,1) -- (0,0,1) -- (1,0,1) -- (1,1,1);
            \draw[blue] (1,0,1) -- (1,0,0) -- (1,1,0);
            \path (0.35,0.35,0.5) node[scale=1,black] {$\bS$};
        \end{scope}
        \begin{scope}[shift={(1,1,2)}]
            \fill[blue!10] (0,1,0) -- (1,1,0) -- (1,1,1) -- (0,1,1) -- cycle;
            \shade[bottom color=blue!60, top color=blue!30] (1,0,0) -- (1,1,0) -- (1,1,1) -- (1,0,1) -- cycle;
            \shade[bottom color=blue!60, top color=blue!30] (0,0,1) -- (1,0,1) -- (1,1,1) -- (0,1,1) -- cycle;
            \draw[blue!70, dashed] (1,0,0) -- (0,0,0) -- (0,1,0);
            \draw[blue!70, dashed] (0,0,0) -- (0,0,1);
            \draw[blue] (0,1,1) -- (0,0,1) -- (1,0,1) -- (1,1,1);
            \draw[blue] (1,0,1) -- (1,0,0) -- (1,1,0);
            \path (0.35,0.35,0.5) node[scale=1,black] {$\bS$};
        \end{scope}
        \begin{scope}[shift={(1,2,1)}]
            \fill[blue!10] (0,1,0) -- (1,1,0) -- (1,1,1) -- (0,1,1) -- cycle;
            \shade[bottom color=blue!60, top color=blue!30] (1,0,0) -- (1,1,0) -- (1,1,1) -- (1,0,1) -- cycle;
            \shade[bottom color=blue!60, top color=blue!30] (0,0,1) -- (1,0,1) -- (1,1,1) -- (0,1,1) -- cycle;
            \draw[blue!70, dashed] (1,0,0) -- (0,0,0) -- (0,1,0);
            \draw[blue!70, dashed] (0,0,0) -- (0,0,1);
            \draw[blue] (0,1,1) -- (0,0,1) -- (1,0,1) -- (1,1,1);
            \draw[blue] (1,0,1) -- (1,0,0) -- (1,1,0);
            \path (0.35,0.35,0.5) node[scale=1,black] {$\bS$};
        \end{scope}
        \begin{scope}[shift={(2,0,1)}]
            \fill[blue!10] (0,1,0) -- (1,1,0) -- (1,1,1) -- (0,1,1) -- cycle;
            \shade[bottom color=blue!60, top color=blue!30] (1,0,0) -- (1,1,0) -- (1,1,1) -- (1,0,1) -- cycle;
            \shade[bottom color=blue!60, top color=blue!30] (0,0,1) -- (1,0,1) -- (1,1,1) -- (0,1,1) -- cycle;
            \draw[blue!70, dashed] (1,0,0) -- (0,0,0) -- (0,1,0);
            \draw[blue!70, dashed] (0,0,0) -- (0,0,1);
            \draw[blue] (0,1,1) -- (0,0,1) -- (1,0,1) -- (1,1,1);
            \draw[blue] (1,0,1) -- (1,0,0) -- (1,1,0);
            \path (0.35,0.35,0.5) node[scale=1,black] {$\bS$};
        \end{scope}
        \draw[draw=blue!60] (3,0,0) -- (3,3,0) -- (3,3,3) -- (3,0,3) -- cycle;
        \draw[draw=blue!60] (0,3,0) -- (0,3,3) -- (0,0,3);
        \draw[draw=blue!60] (0,3,3) -- (3,3,3);
        \path (-1.5,1.5,2) node {$\bG = $};
    \end{tikzpicture}
    \caption{\emph{Structure of the core tensor $\bG$ for the case $d=3$ resulting from an application of the anisotropic diffusion operator.}}
    \label{fig:structured_tens}
\end{figure}
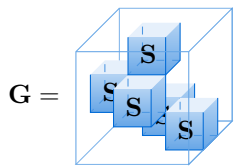

The corresponding linear operator $\mcA$ acting on $\bX \in \R^{n_1 \times \cdots \times n_d}$ can be represented as a TT operator of rank 3, with the cores given by
\[
    A_1(i_1,j_1) =  \begin{bmatrix} 
                    L_{1}(i_1,j_1) & 2 \alpha B_1(i_1,j_1) & I_{n_1}(i_1,j_1) 
                \end{bmatrix},
                \quad A_d(i_d,j_d) = \begin{bmatrix} 
                    I_{n_d}(i_d,j_d) \\ B_d(i_d,j_d) \\ L_d(i_d,j_d) 
                \end{bmatrix},
\]
and
\[
    A_\mu(i_\mu,j_\mu) = \begin{bmatrix} 
                            I_{n_\mu}(i_\mu,j_\mu) & 0 & 0 \\ 
                            B_\mu(i_\mu,j_\mu) & 0 & 0  \\ 
                            L_\mu(i_\mu,j_\mu) & 2\alpha B_\mu(i_\mu,j_\mu) & I_{n_\mu}(i_\mu,j_\mu) 
                        \end{bmatrix}, \quad \mu = 2,\dots,d-1.
\]
In the Tucker format, this operator is also of rank $3$. Given a tensor $\bX$ in the representation \eqref{eq:tuckerformat}, the result $\bY = \mcA \bX$ is explicitly given by $\bY = \bG \times_1 V_1 \times_2 \cdots \times_d V_d$ with
\[
    V_\mu = \begin{bmatrix}
        U_\mu & L_\mu U_\mu & B_\mu U_\mu 
        \end{bmatrix} \in \R^{n \times 3r_\mu}
\]
and core tensor $\bG \in \R^{3r_1 \times \cdots \times 3r_d}$ which has a block structure shown in Figure \ref{fig:structured_tens} for the case $d=3$.

The rank of $\mcA$ increases linearly with the band width of the diffusion matrix $D$. For example, a pentadiagonal structure would yield an operator of rank $4$. See also~\cite{Kazeev:2013} for more general bounds in terms of certain properties of $D$.

\subsection{Results for the Tucker format}
\label{sec:experiments_tucker}

\begin{figure}[h]
    \centering
    \includegraphics[width=0.48\textwidth]{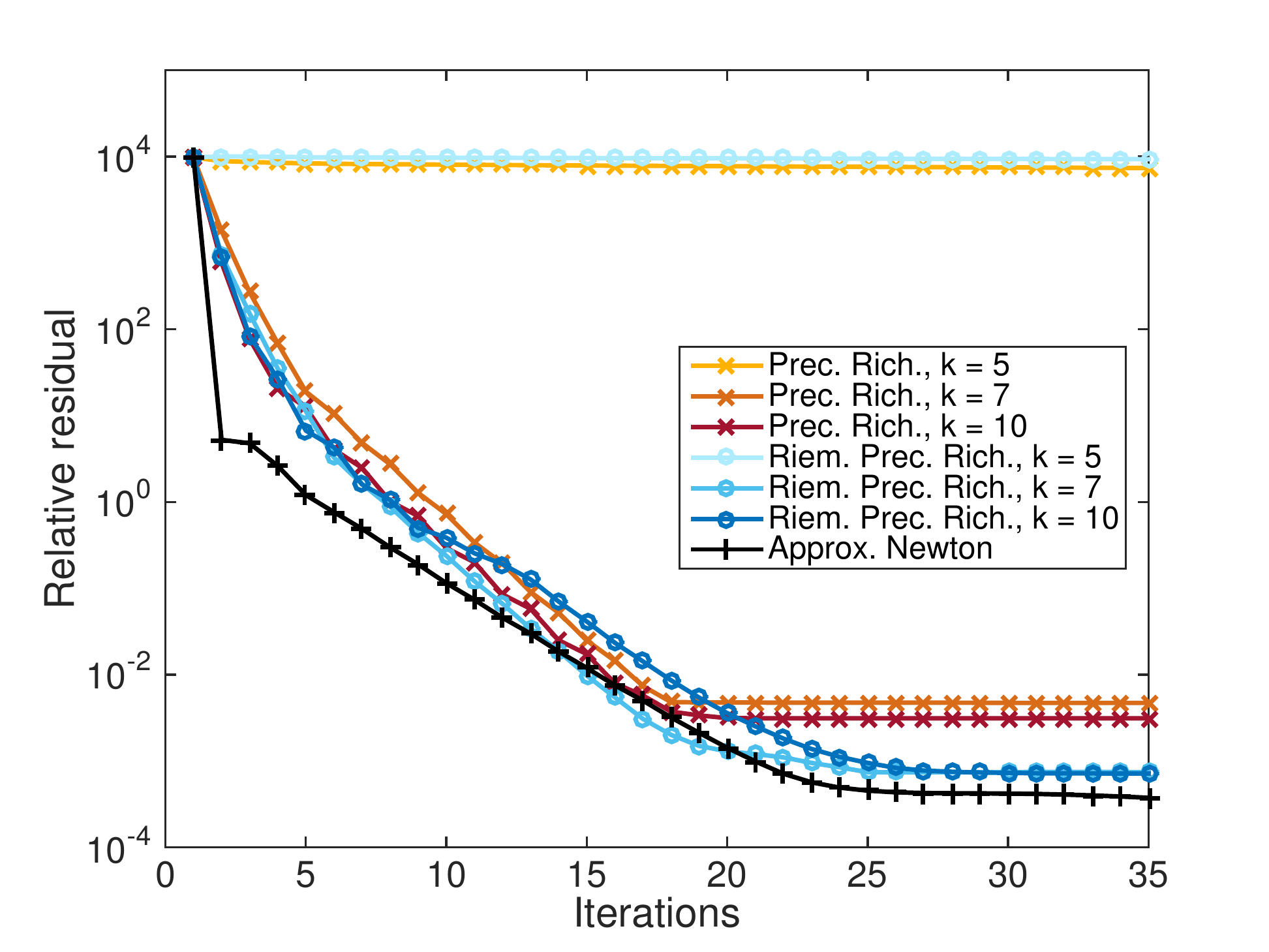}
    \includegraphics[width=0.48\textwidth]{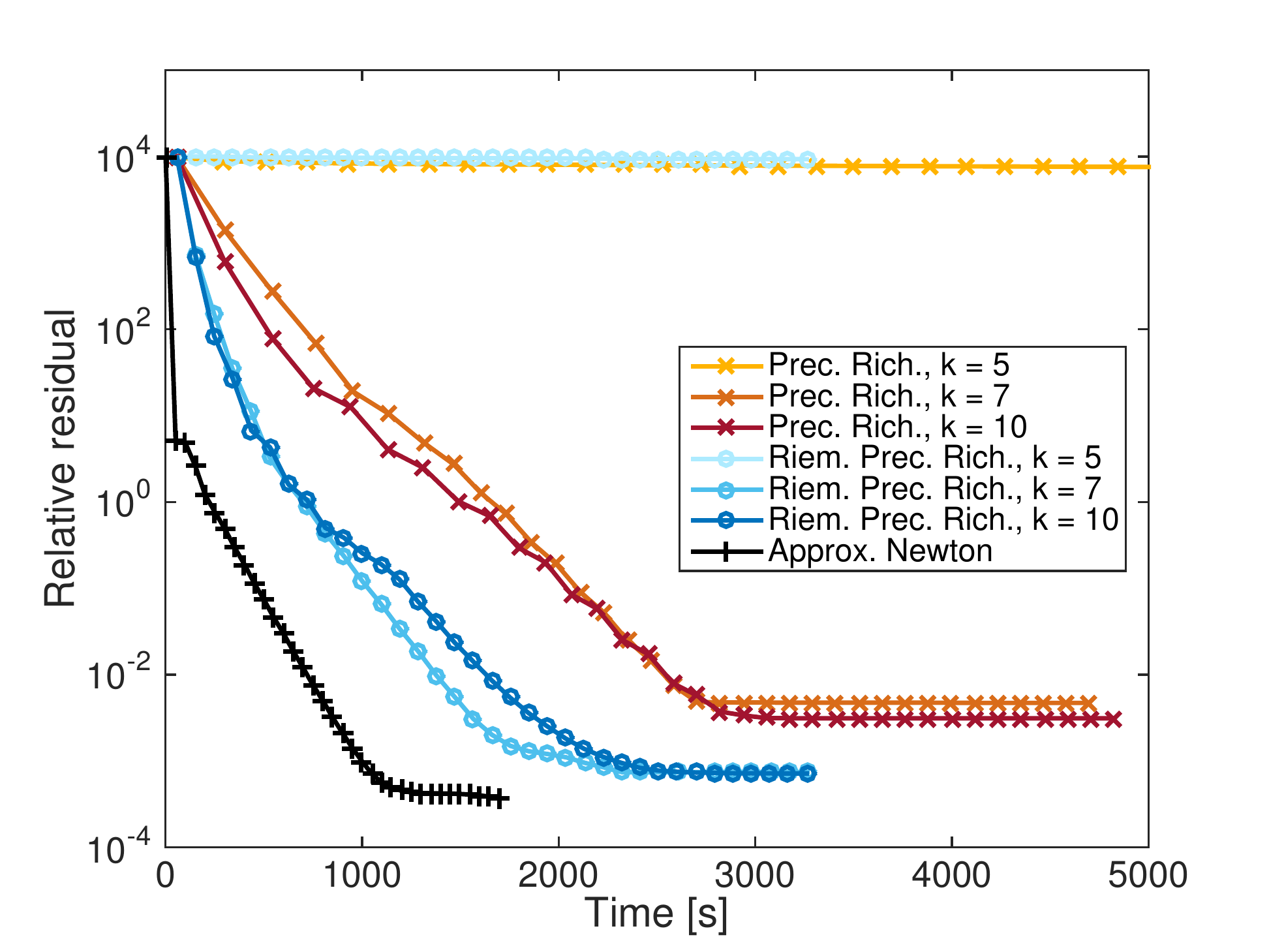}
    \caption{Newton potential with $d = 3$. \emph{Comparison of truncated preconditioned Richardson,
    truncated Riemannian preconditioned Richardson, and the approximate Newton scheme when
    applied to the Newton potential in the Tucker format. For the Richardson iterations, exponential sum approximations with
    $k \in \{5,7,10\}$ terms are compared. \textbf{Left:} Relative residual as a function of
    iterations. \textbf{Right:} Relative residual as a function of time} }
    \label{fig:tucker_newton}
\end{figure}

\begin{figure}[h]
    \centering
    \includegraphics[width=0.48\textwidth]{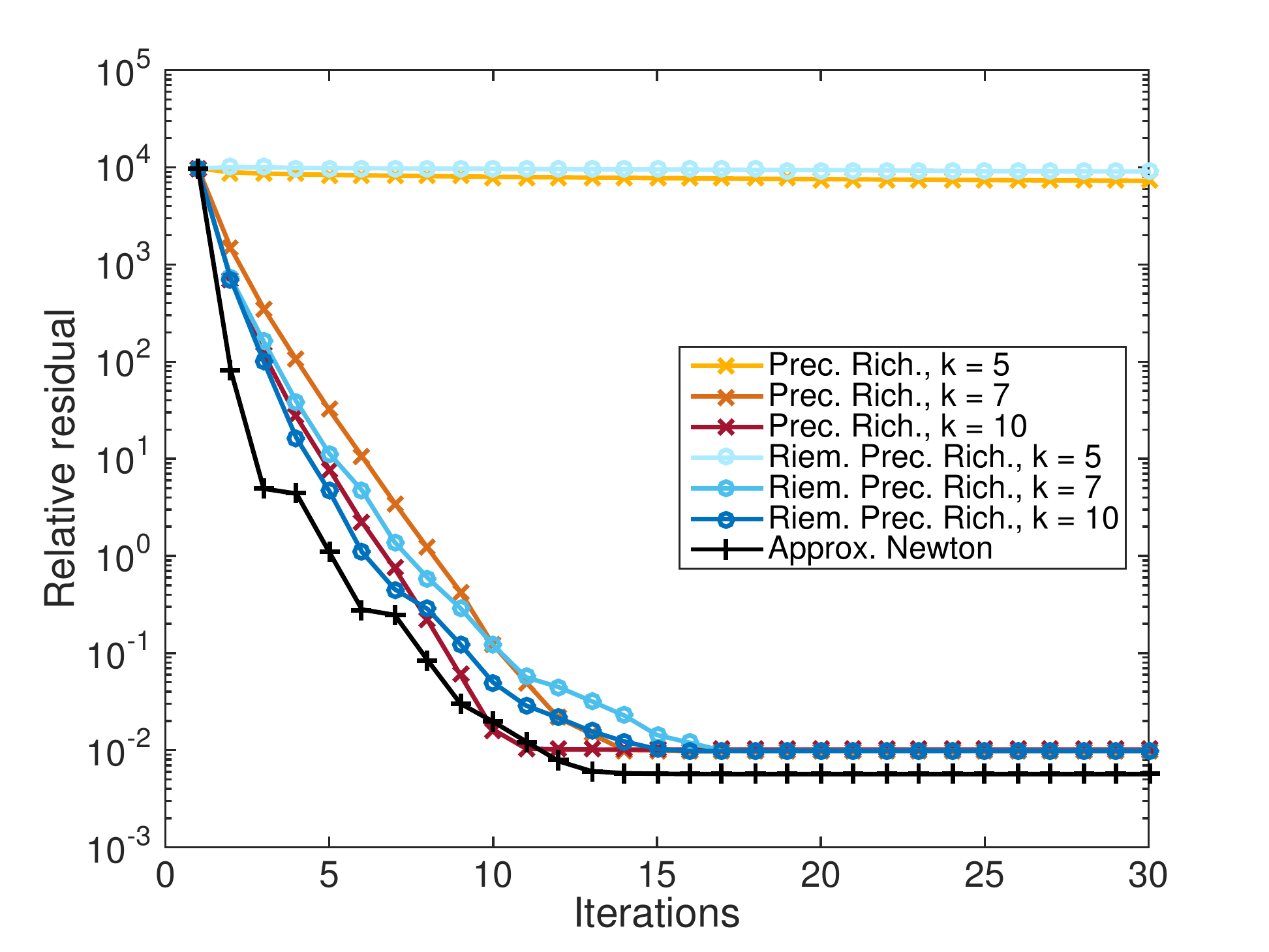}
    \includegraphics[width=0.48\textwidth]{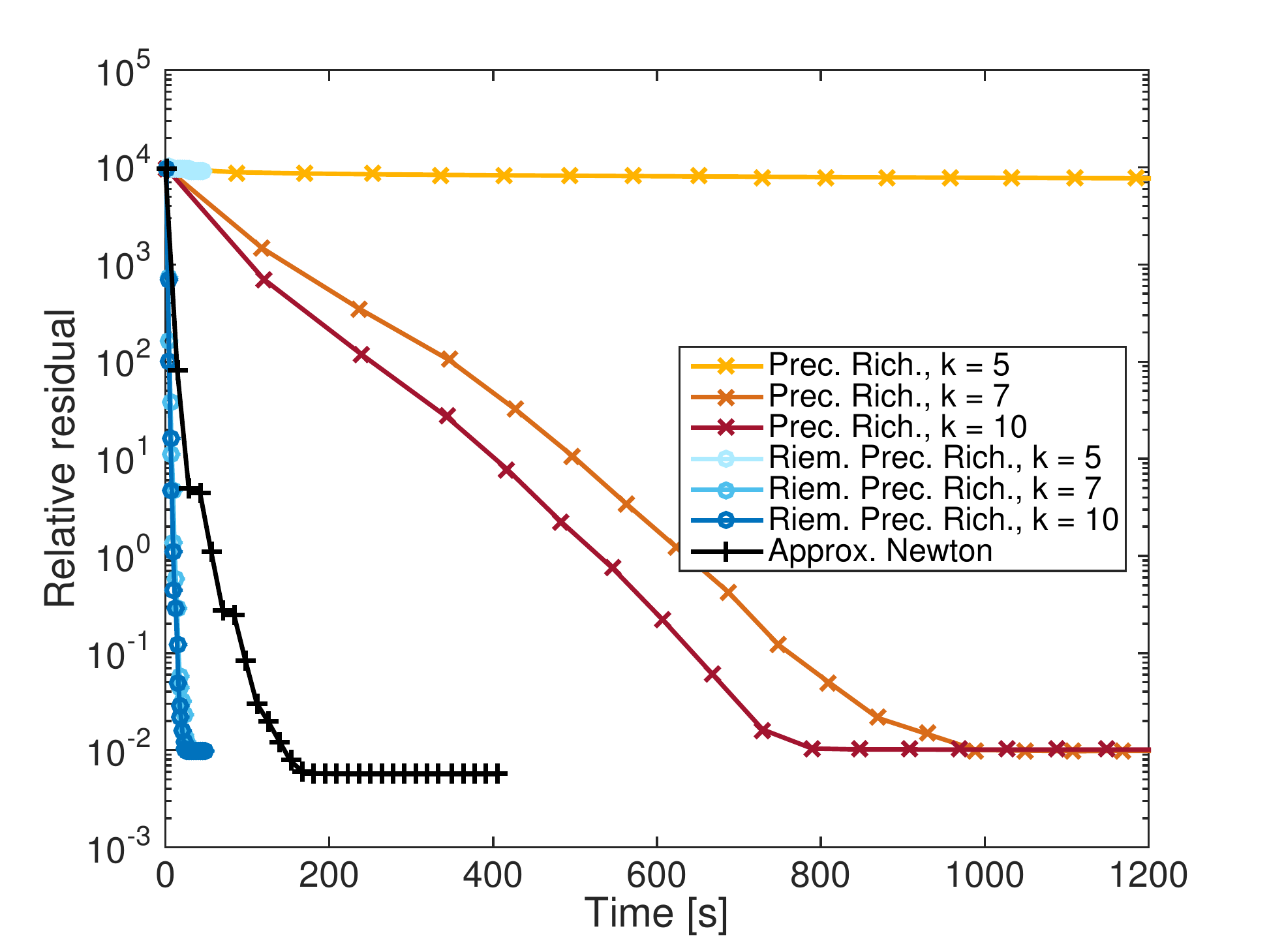}
    \caption{Anisotropic diffusion with $d = 3$. \emph{Comparison of truncated Preconditioned Richardson,
    truncated Riemannian preconditioned Richardson, and the approximate Newton scheme when
    applied to the Newton potential in the Tucker format. For the Richardson iterations, exponential sum approximations with
    $k \in \{5,7,10\}$ terms are compared. \textbf{Left:} Relative residual as a function of
    iterations. \textbf{Right:} Relative residual as a function of time} }
    \label{fig:tucker_anisotropic}
\end{figure}

For tensors represented in the Tucker format we want to investigate the convergence of the truncated
preconditioned Richardson \eqref{eq:trunc_prec_Richardson} and its Riemannian variant
\eqref{eq:trunc_prec_projected_Richardson}, and compare them to the approximate Newton scheme 
discussed in \S\ref{sec:newton_tucker}.
Figure~\ref{fig:tucker_newton} displays the obtained results for the first test case, the Newton potential, 
where we set $d=3$, $n = 100$, and used multilinear ranks $r=15$.
Figure \ref{fig:tucker_anisotropic} displays the results for the second test case, the anisotropic diffusion operator with $\alpha = \frac{1}{4}$, using the same settings.
In both cases, the right hand side is given by a random rank-1 Tucker tensor. 
To create a full space preconditioner for both Richardson approaches, we approximate the inverse
Laplacian by an exponential sum of $k \in \{5,7,10\}$ terms. It can be clearly  seen that the
quality of the preconditioner has a strong influence on the convergence. For
$k=5$, convergence is extremely slow. Increasing $k$ yields a drastic improvement on the convergence.

With an accurate preconditioner, the truncated Richardson scheme converges fast
with regard to the number of iterations, but suffers from very long computation
times due to the exceedingly high intermediate ranks. In comparison, the
Riemannian Richardson scheme yields similar convergence speed, but with
significantly reduced computation time due to the additional projection into
the tangent space. The biggest saving in computational effort comes from
relation~\eqref{eq:proj_commute_prec} which allows us to avoid having to form
the preconditioned residual $\mcP^{-1} (\bF - \mcA \bX_k)$ explicitly, a
quantity with very high rank.  Note that for both Richardson approaches, it is
necessary to round the Euclidean gradient to lower rank using a tolerance of,
say, $10^{-5}$ before applying the preconditioner to avoid excessive
intermediate ranks.

The approximate Newton scheme converges equally well as the best Richardson
approaches with regard to the number of iterations and does not require setting up a preconditioner.
For the first test case, it only needs about half of the time as the best
Richardson approach. For the second test case, 
it is significantly slower than Riemannian preconditioned Richardson. Since this
operator is of lower rank than the Newton potential, the additional complexity
of constructing the approximate Hessian does not pay off in this case.

\paragraph{Quadratic convergence.}

\begin{figure}[h]
    \centering
    \includegraphics[width=0.48\textwidth]{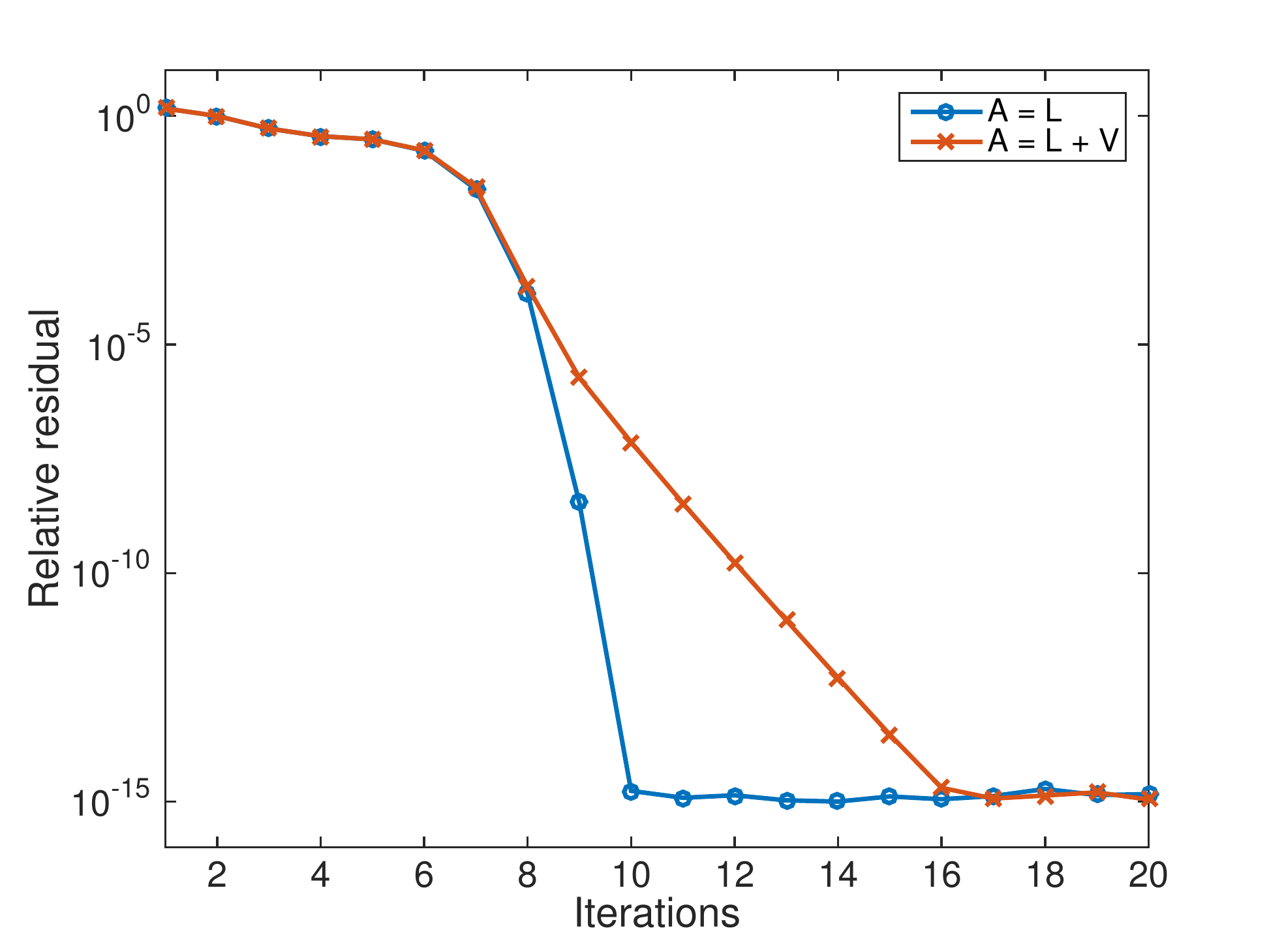}
    \caption{\emph{Convergence of the approximate Newton method for the zero-residual case when applied to a pure Laplace
    operator $L$ and to the anisotropic diffusion operator $L + V$.}}
    \label{fig:quad_conv}
\end{figure}

In Figure \ref{fig:quad_conv} we investigate the convergence of the approximate
Newton scheme when applied to a pure Laplace operator, $A = L$, and to the
anisotropic diffusion operator $A = L + V$.  In order to have an exact solution
of known rank $r = 4$, we construct the right hand side by applying $A$ to a
random rank $4$ tensor. For the dimension and tensor size we have chosen $d=3$
and $n = 200$, respectively.  By construction, the exact solution lies on the
manifold. Hence, if the approximate Newton method converges to this solution,
we have zero residual and our Gauss--Newton approximation of~\eqref{eq:Hessian}
is an exact second-order model despite only containing the $A$ term. In other
words, we expect quadratic convergence when $A = L$ but only linear when $A =
L+V$ since our approximate Newton method~\eqref{eq:GN_manifold} only solves
with $L$. This is indeed confirmed in Figure~\ref{fig:quad_conv}. 

\subsection{Results for the TT format}

In the TT format, we compare the convergence of our approximate Newton scheme (with the overlapping block-Jacobi
preconditioner described in \S\ref{sec:overlapping}) to a standard 
approach, the alternating linear scheme (ALS).

\begin{figure}[h]
    \includegraphics[width=0.48\textwidth]{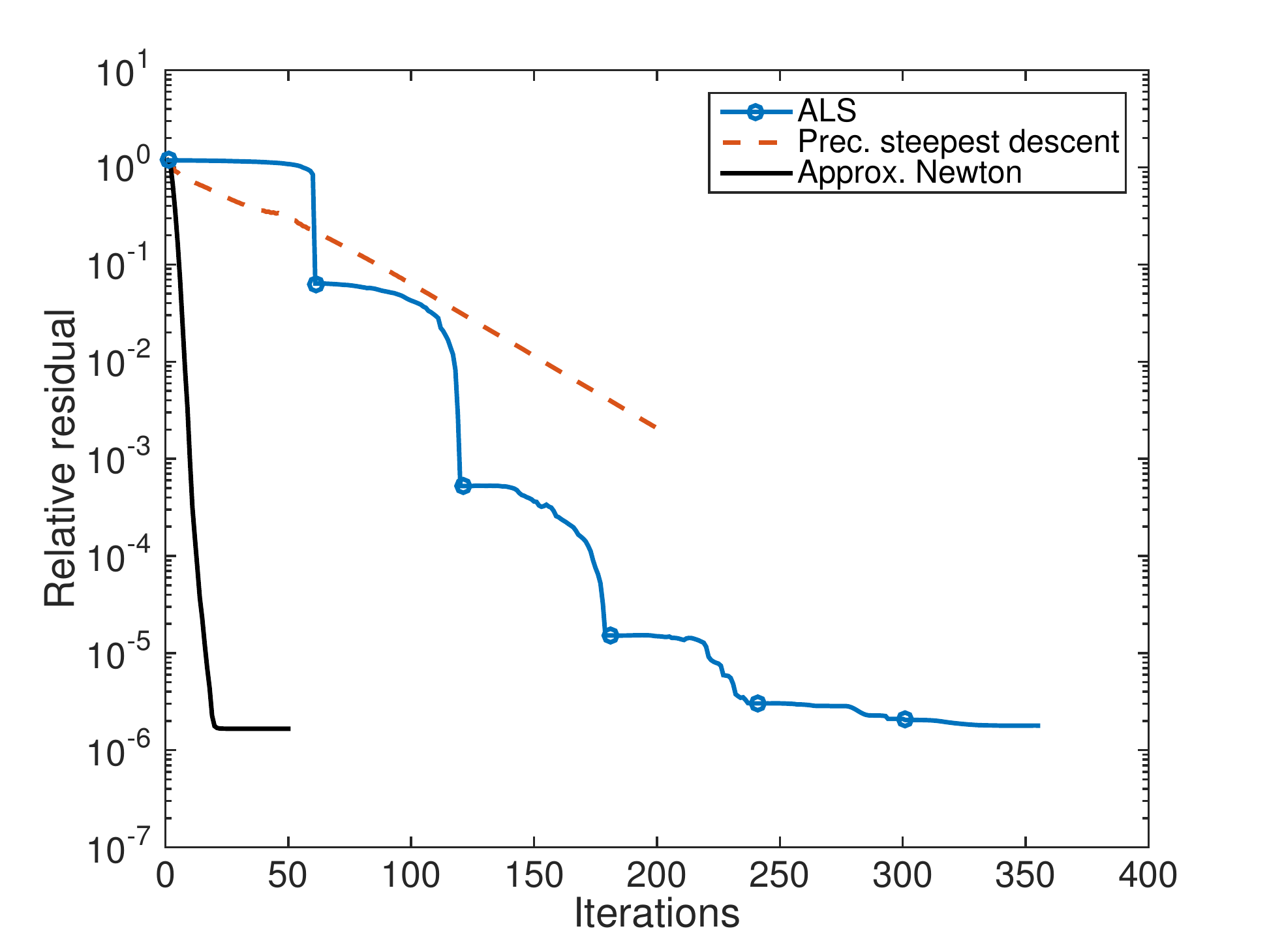}
    \includegraphics[width=0.48\textwidth]{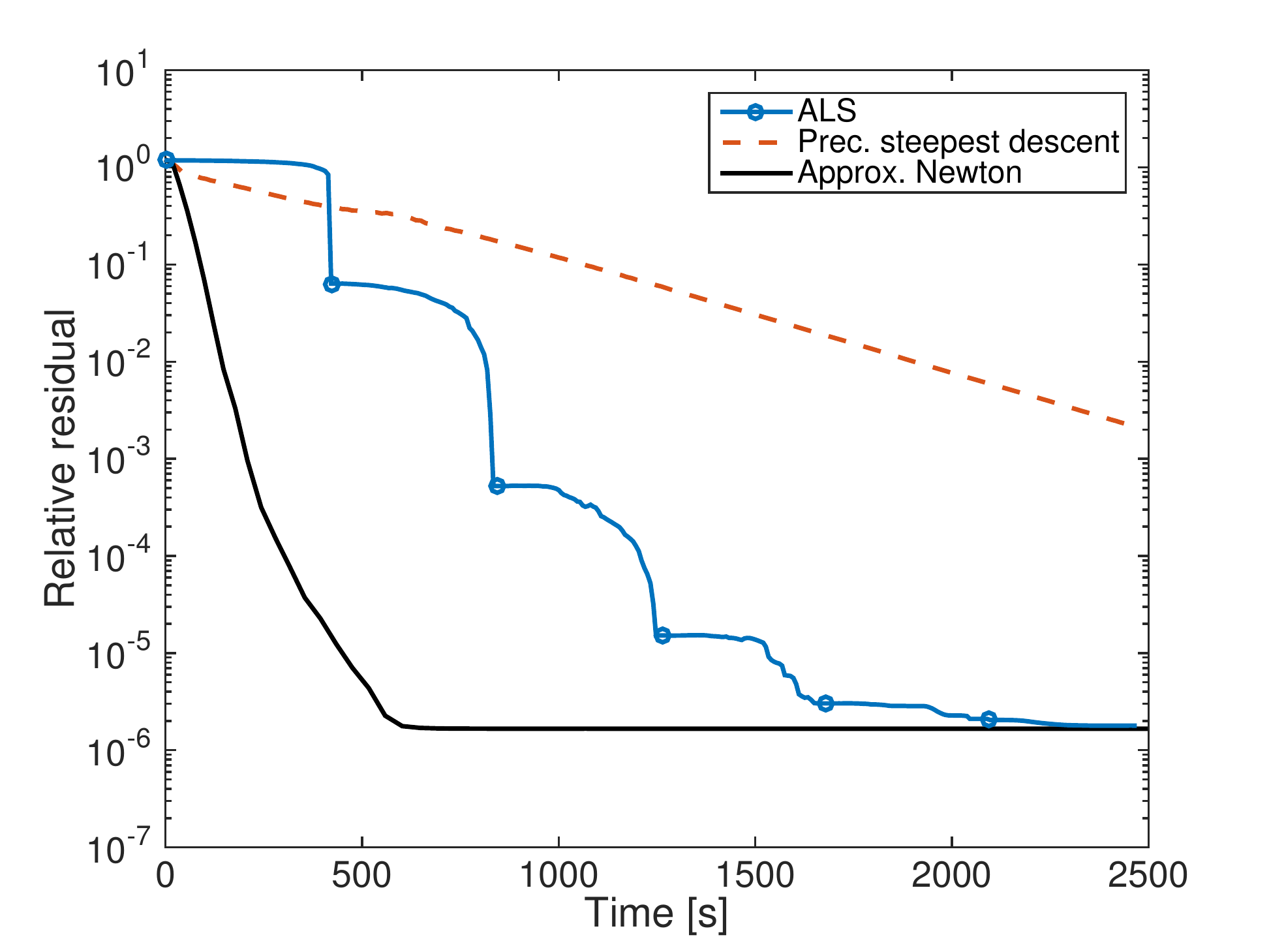}
        \caption{Newton potential with $d = 60$. \emph{Convergence of ALS compared to
                preconditioned steepest descent with overlapping block-Jacobi
                as preconditioner and the approximate Newton scheme.
                \textbf{Left:} Relative residual as function of iterations.
                \textbf{Right:} Relative residual as function of time. }}
    \label{fig:TT_newton}
\end{figure}

We have chosen $d = 60$, $n = 100$, and
a random rank-one right hand sides of norm one. In the
first test case, the Newton potential, we have chosen TT ranks $r =
10$ for the approximate solution. The corresponding convergence curves are shown in Figure
\ref{fig:TT_newton}. We observe that the approximate Newton scheme needs
significantly less time to converge than the ALS scheme. As a reference, we
have also included a steepest descent method using the overlapping block-Jacobi scheme
 directly as a preconditioner for every gradient step instead of using it to
 solve the approximate Newton equation \eqref{eq:GN_TT}. The additional effort of solving the Newton equation approximately clearly pays off.

\begin{figure}[h]
    \includegraphics[width=0.48\textwidth]{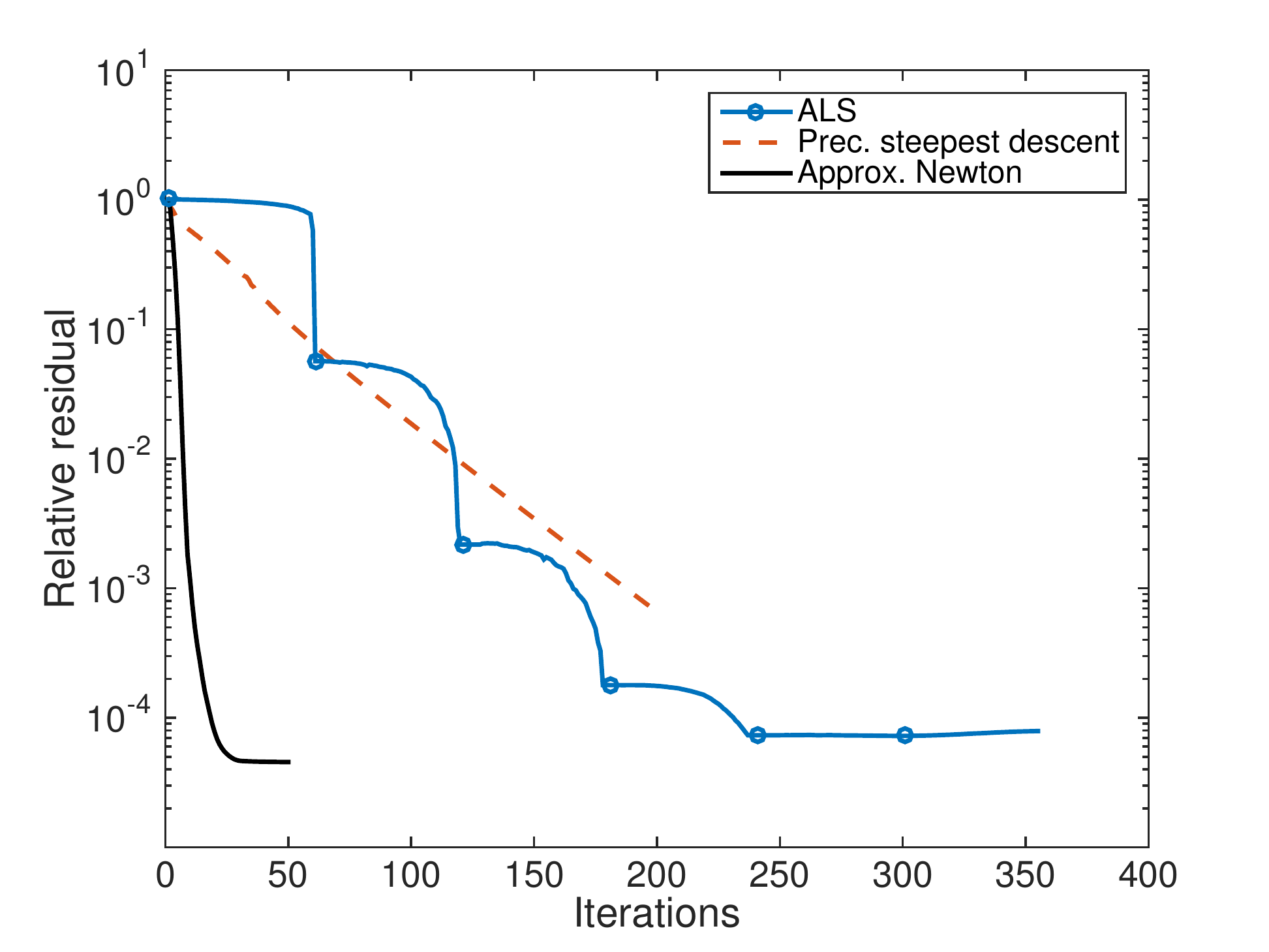}
    \includegraphics[width=0.48\textwidth]{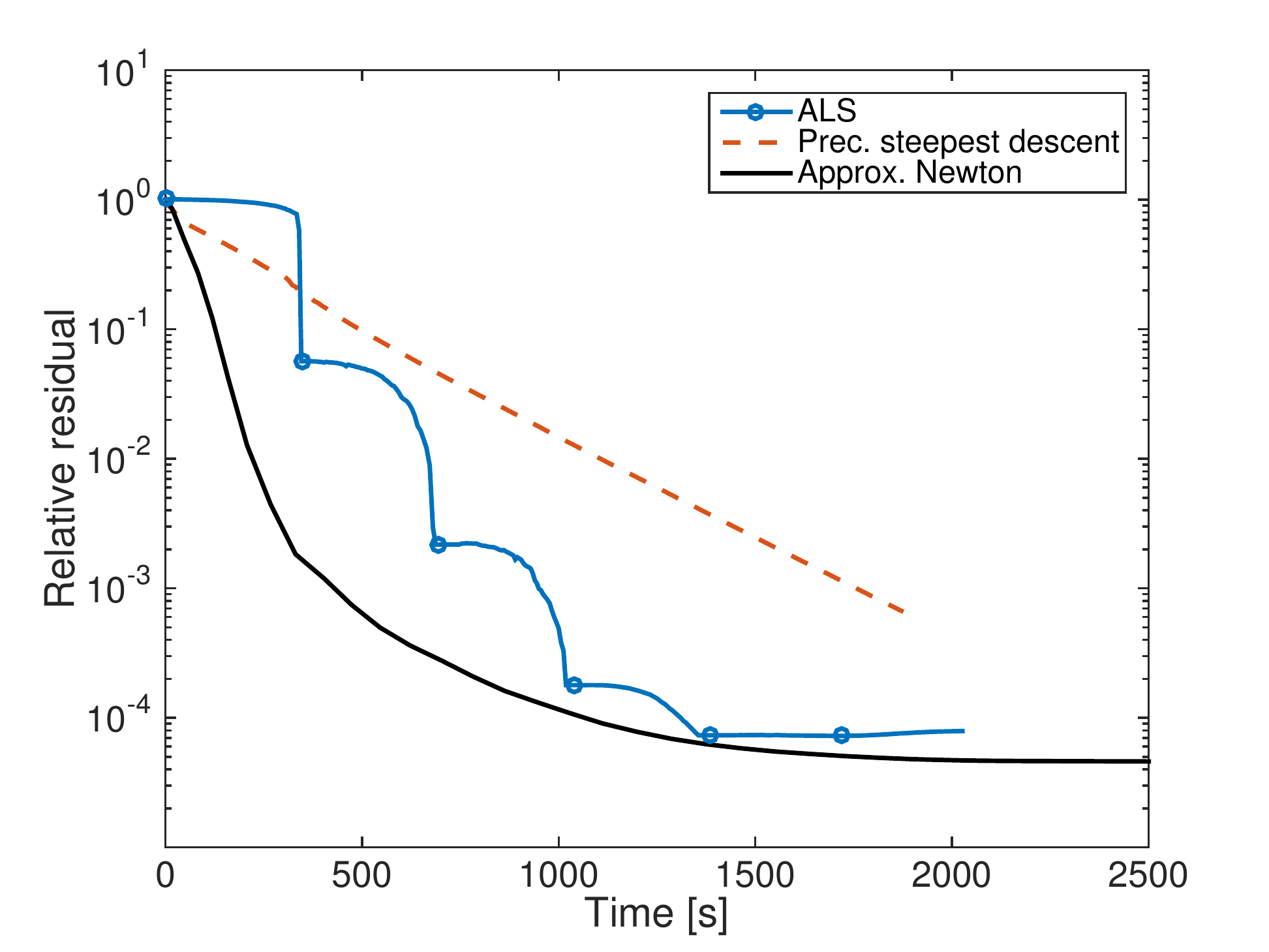}
        \caption{Anisotropic diffusion with $d = 60$. \emph{Convergence of ALS compared to
                preconditioned steepest descent with overlapping block-Jacobi
                as preconditioner and the approximate Newton scheme.
                \textbf{Left:} Relative residual as function of iterations.
                \textbf{Right:} Relative residual as function of time. }}
    \label{fig:TT_diffusion}
\end{figure}

In Figure \ref{fig:TT_diffusion}, we show results for the anisotropic diffusion
case. To obtain a good accuracy of  the solution, we have to choose a
relatively high rank of $r = 25$ in this case. Here, the approximate Newton
scheme is still faster, especially at the beginning of the iteration, but the
final time needed to reach a residual of $10^{-4}$ is similar to ALS.

Note that in Figures \ref{fig:TT_newton} and \ref{fig:TT_diffusion} the plots
with regard to the number of iterations are to be read with care due to the
different natures of the algorithms. One ALS iteration corresponds to the 
optimization of one core. In the plots, the beginning of each half-sweep of ALS is denoted by a circle.
To assessment the performance of both schemes as fairly as possible, we
have taken considerable care to provide the same level of optimization to the
implementations of both the ALS and the approximate Newton scheme.

\paragraph{Mesh-dependence of the preconditioner.} 
\begin{figure}[h]
    \centering
    \includegraphics[width=0.48\textwidth]{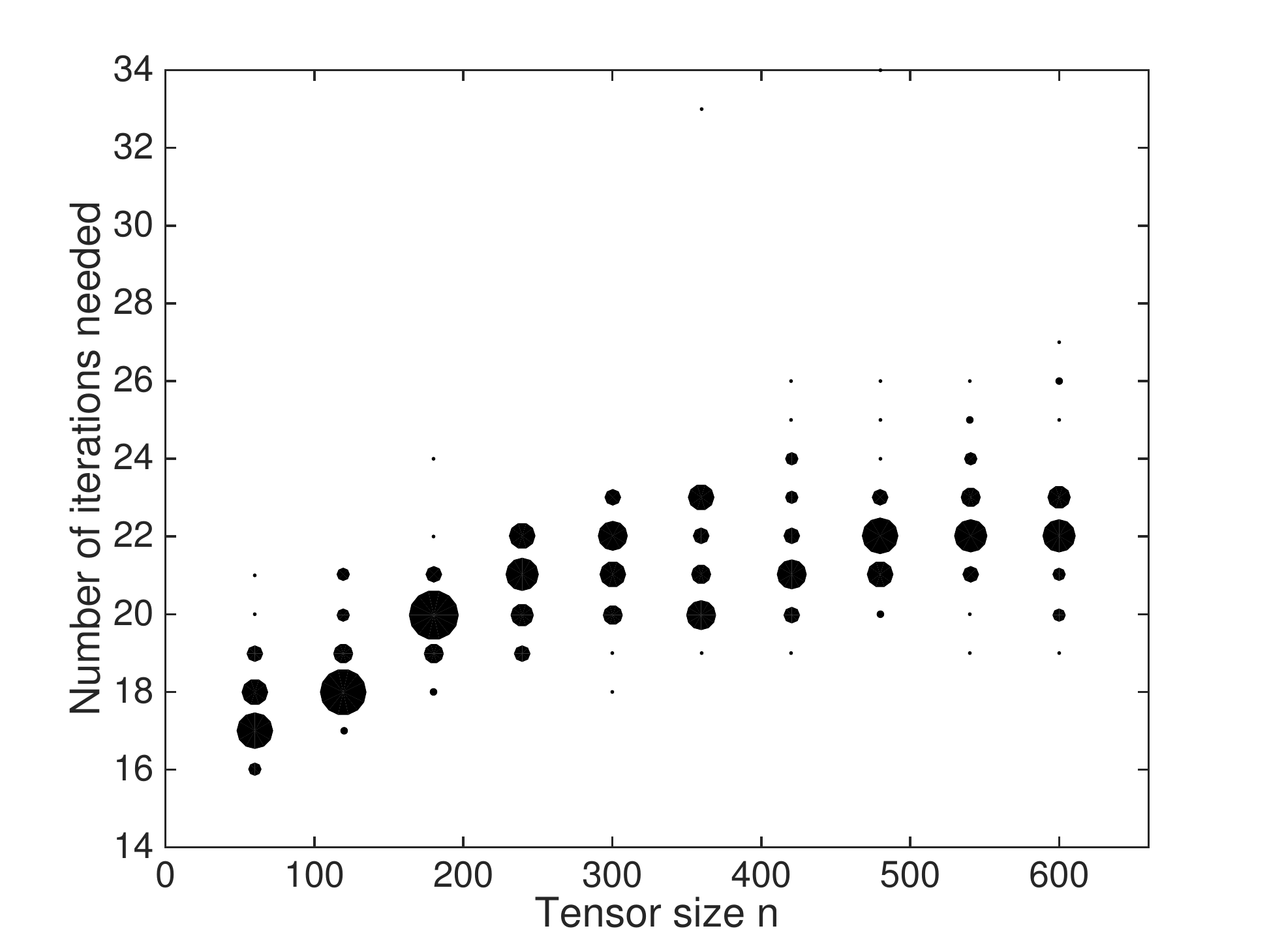}
    \includegraphics[width=0.48\textwidth]{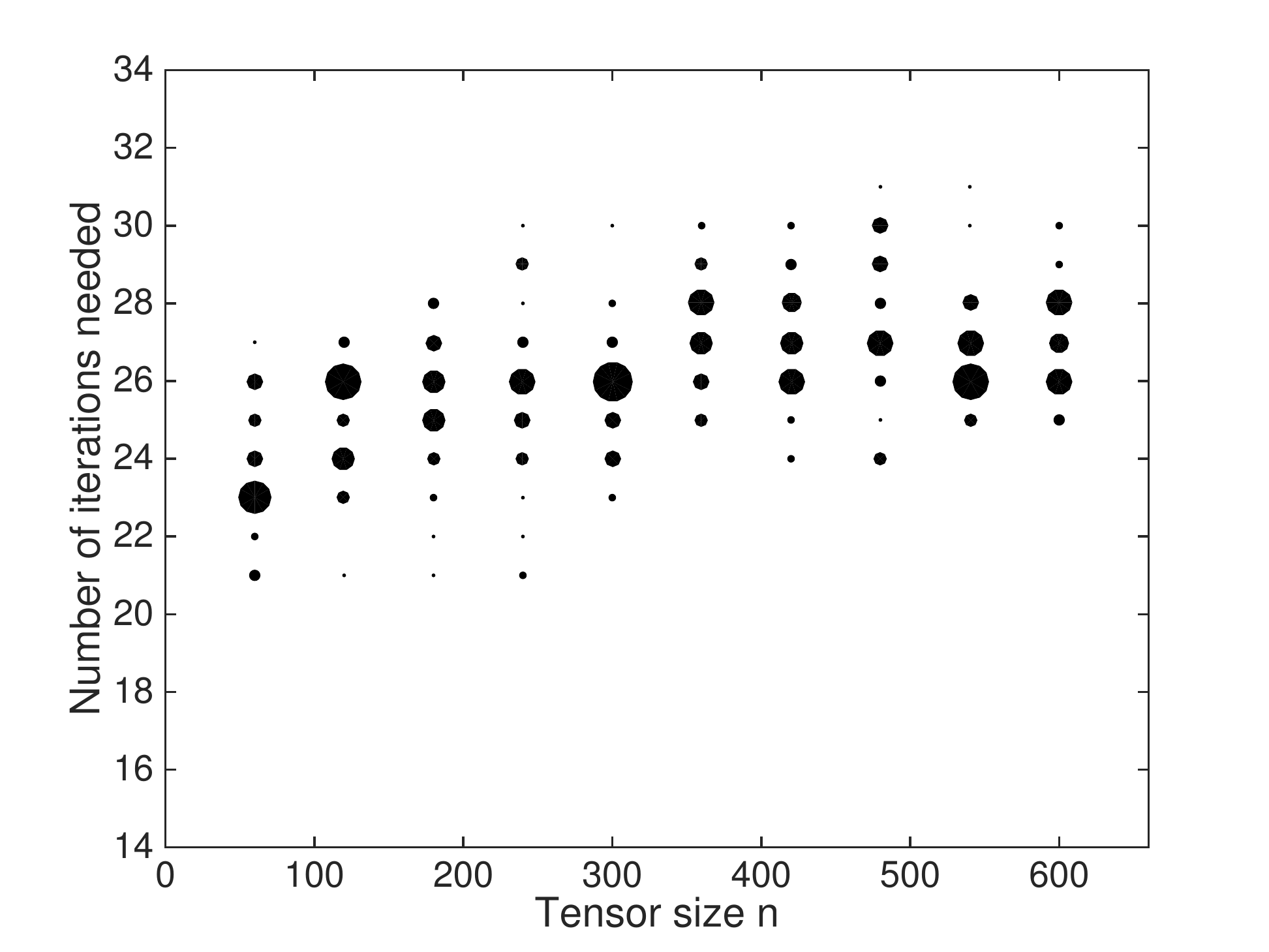}
    \caption{\emph{Number of iterations that the proposed approximate Newton
            scheme needs to reach a relative residual of $10^{-6}$ for
            different mesh widths $h = 1/n$. The solution has dimension $d=10$
            and rank $r = 3$.  We perform 30 runs for each size. The radii of
            the circles corresponds to the number of runs achieving this number
            of iterations. \textbf{Left:} Dimension $d=10$. \textbf{Right:}
            Dimension $d=30$.}}
    \label{fig:meshindep}
\end{figure}

To investigate how the
performance of the preconditioner depends on the mesh width of the
discretization, we look again at the anisotropic diffusion operator  and
measure the convergence as the mesh width $h$ and therefore the tensor size $n
\in \{60,120,180,240,360,420,480,540,600\}$ changes by one order of magnitude.
As in the test for quadratic convergence, we construct the right hand side by
applying $A$ to a random rank $3$ tensor. For the dimension and tensor size we
have chosen $d=3$ and $n = 200$, respectively. 

To measure the convergence, we take the number of iterations needed to converge
to a relative residual of $10^{-6}$. For each tensor size, we perform 30 runs
with random starting guesses of rank $r = 3$. The result is shown in Figure
\ref{fig:meshindep}, where circles are drawn for each combination of size
$n$ and number of iterations needed. The radius of each circle denotes how many
runs have achieved a residual of $10^{-6}$ for this number of iterations. 

On the left plot of \ref{fig:meshindep} we see the results of dimension $d=10$,
whereas on the right plot we have $d=30$. We see that the number of iterations
needed to converge changes only mildly as the mesh width varies over one order
of magnitude. In addition, the dependence on $d$ is also not very large.

\subsection{Rank-adaptivity}
\begin{figure}[h]
    \centering
    \includegraphics[width=0.55\textwidth]{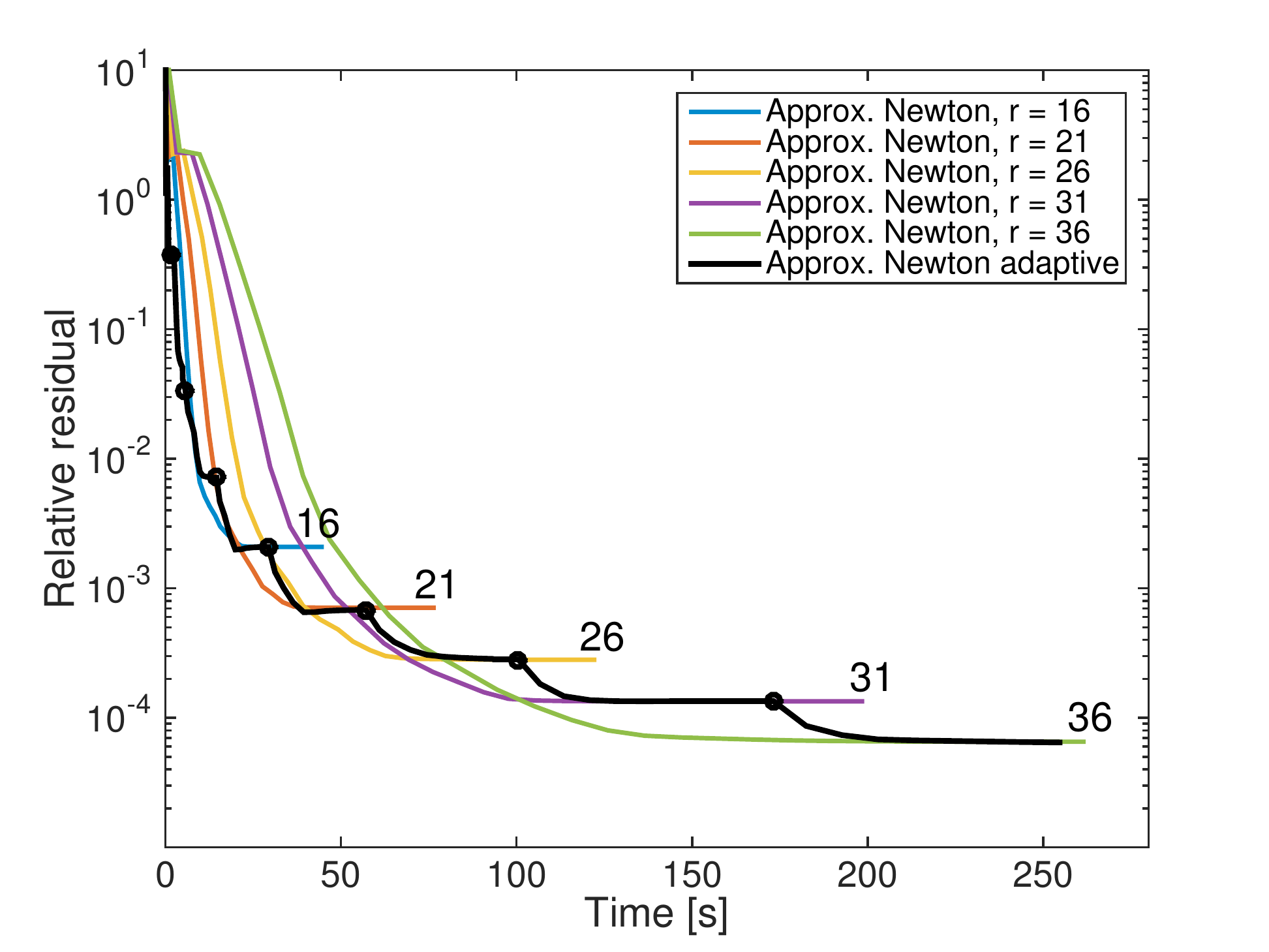}
    \caption{\emph{Rank-adaptivity for approximate Newton applied to the
            anisotropic diffusion equation with $n=100$, $d=10$. Starting from
            rank 1, the rank is increased by 5 after 10 iterations per rank.
            Each rank increase is denoted by a black circle. The other curves
            show the convergence when running approximate Newton directly with
            the target rank.}}
    \label{fig:adaptive}
\end{figure}

Note that in many applications, rank-adaptivity of the algorithm is a desired
property. For the Richardson approach, this would result in replacing the
fixed-rank truncation with a tolerance-based rounding procedure. In the
alternating optimization, this would lead to the DMRG or AMEn algorithms. 
In the framework of Riemannian optimization, rank-adaptivity can be introduced
by successive runs of increasing rank, using the previous solution as a warm
start for the next rank. For a recent discussion of this approach, see
\cite{Uschmajew_V_2015}. A basic example of introducing rank-adaptivity to the
approximate Newton scheme is shown in Figure \ref{fig:adaptive}. Starting from
ranks $r^{(0)} = 1$, we run the approximate Newton scheme for 10 iterations and use this
result to warm start the algorithm with ranks $r^{(i)} = r^{(i-1)} + 5$. At each
rank, we perform 10 iterations of the approximate Newton scheme. The result is
compared to the convergence of approximate Newton when starting directly with
the target rank $r^{(i)}$. We see that the obtained relative residuals match for each of the ranks $r^{(i)}$. Although the adaptive rank scheme is slower for a
desired target rank due to the additional intermediate steps, it offers more
flexibility when we want to instead prescribe a desired accuracy. For a relative
residual of $10^{-3}$, the adaptive scheme needs about half the time than using
the (too large) rank $r=36$. 

Note that in the case of tensor
completion, rank adaptivity becomes a crucial ingredient to avoid overfitting
and to steer the algorithm into the right direction, see e.g.
\cite{Vandereycken2012,Kressner2013b,Tan2014,Uschmajew_V_2015,Steinlechner2015}. For difficult completion problems,
careful core-by-core rank increases become necessary. Here, for linear systems,
such a core-by-core strategy does not seem to be necessary, as the algorithms will converge even if we
directly optimize using rank $r=36$. This is likely due to the preconditioner which acts globally over all cores.

\FloatBarrier

\section{Conclusions}

We have investigated different ways of introducing preconditioning into
Riemannian gradient descent. As a simple but effective approach, we have seen
the Riemannian truncated preconditioned Richardson scheme. Another approach used second-order
information by means of approximating the Riemannian Hessian. In the
Tucker case, the resulting approximate Newton equation could be solved
efficiently in closed form, whereas in the TT case, we have shown that this equation can be
solved iteratively in a very efficient way using PCG with an overlapping
block-Jacobi preconditioner. The numerical experiments show favorable
performance of the proposed algorithms when compared to standard non-Riemannian approaches, such as truncated preconditioned Richardson and ALS. The advantages of the approximate Newton scheme become especially pronounced in cases when the linear operator is expensive to apply, e.g., the Newton potential.
 
\bibliographystyle{plain}
\bibliography{linearsystem}
\end{document}